\documentclass[12pt]{article}
\usepackage[final]{epsfig}
\usepackage{graphics}
\usepackage{amsmath,amssymb,amsthm,amscd}
\usepackage{latexsym}
\usepackage{epsfig}
\usepackage{amsfonts}
\usepackage{latexsym}
\usepackage{graphicx}
\usepackage{epstopdf}
\usepackage{multicol,multirow}
\def\diag{\text{diag} \thinspace}

\def\proof{\paragraph{Proof.}}
\def\prooflem{\paragraph{Proof of lemma.}}
\def\proofthm{\paragraph{Proof of theorem.}}
\def\proofthms{\paragraph{Proof of theorems}}

\def\bull{\noindent $\bullet$ \thinspace}

\newtheorem{lemma}{Lemma}[section]
\newtheorem{proposition}[lemma]{Proposition}
\newtheorem{remark}[lemma]{Remark}
\newtheorem{remark-definition}[lemma]{Remark-Definition}
\newtheorem{example}[lemma]{Example}
\newtheorem{theorem}[lemma]{Theorem}
\newtheorem{definition}[lemma]{Definition}
\newtheorem{corollary}[lemma]{Corollary}

\newtheorem{proposition-conjecture}[lemma]{Proposition-conjecture}
\newtheorem{problem}[lemma]{Problem}

\newcommand{\proofend}{\hfill$\Box$\bigskip}
\newcommand{\C}{{\mathbb C}}

\newcommand{\R}{{\mathbb R}}
\newcommand{\rr}{{\mathbf{r}}}
\newcommand{\Z}{{\mathbb Z}}
\newcommand{\N}{{\mathbb N}}
\newcommand{\RP}{{\mathbb {RP}}}
\newcommand{\CP}{{\mathbb {CP}}}
\newcommand{\PP}{{\mathbb {P}}}
\newcommand{\ep}{\epsilon}

\newcommand{\marginnote}[1]
{
}

\newcounter{bk}

\newcounter{fs}

\title {The geometry of dented pentagram maps}

\author{Boris Khesin 
  and Fedor Soloviev\thanks{Department of Mathematics, University of Toronto, Toronto, ON M5S 2E4, Canada; e-mails: \tt{khesin@math.toronto.edu} and \tt{soloviev@math.toronto.edu}
  }}

\date{}

\begin{document}
\maketitle

\begin{abstract}
We propose a new family of natural generalizations of the pentagram map from 2D
to higher dimensions and prove their integrability on generic twisted and closed polygons.
In dimension $d$ there are $d-1$
such generalizations called dented pentagram maps,
and we describe their geometry, continuous limit,  and Lax representations
with a spectral parameter. We prove algebraic-geometric integrability of the dented pentagram maps
in the 3D case
and compare the dimensions of  invariant tori for the dented maps with those
for the higher pentagram maps constructed with the help of
short diagonal hyperplanes.
When restricted to corrugated polygons, the dented pentagram maps coincide
between themselves and with  the corresponding corrugated pentagram map.
Finally, we prove integrability for a  variety of  pentagram maps for generic 
and partially corrugated polygons in higher dimensions.
\end{abstract}
\tableofcontents

\section*{Introduction} \label{intro}

\begin{figure}[hbtp]
\centering
\includegraphics[width=1.8in]{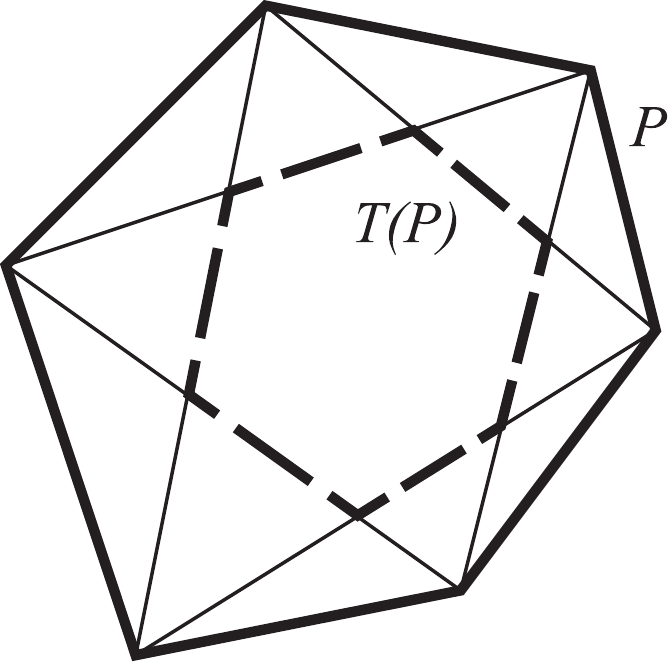}
\caption{\small The image $T(P)$ of a hexagon $P$ under the 2D pentagram map.}
\label{fig:hex}
\end{figure}

The pentagram map was originally defined  in \cite{Schwartz} as a map on plane convex 
polygons considered up to their projective equivalence,
where a new polygon is spanned by the shortest diagonals of the initial one, see Figure \ref{fig:hex}.
This map is the identity for pentagons, it is an involution for hexagons, while 
for polygons with more vertices it was shown to exhibit quasi-periodic behaviour under iterations.
The pentagram  map was extended to the case of twisted polygons and its integrability in 2D was proved in  \cite{OST99}, see also \cite{FS}.

\begin{figure}[hbtp!]
\centering
\includegraphics[width=3.7in]{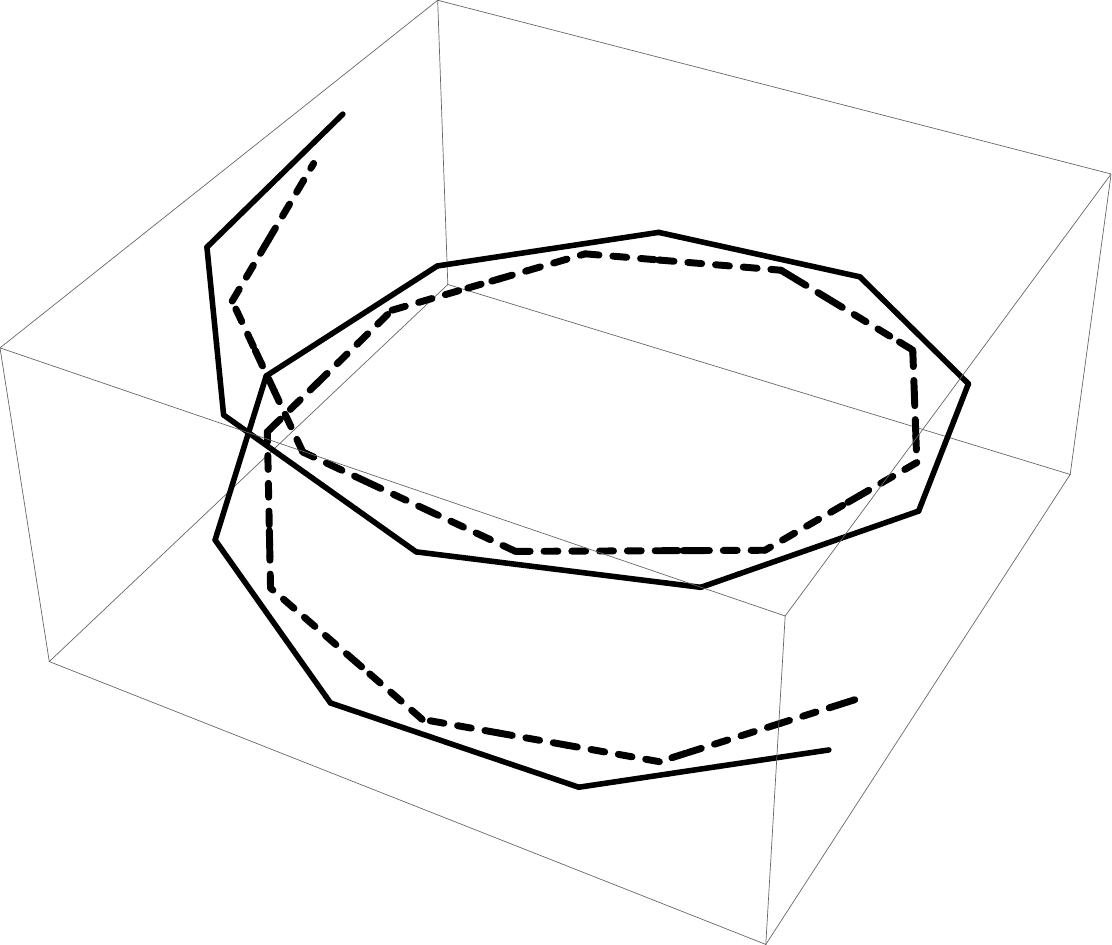}
\caption{\small A space pentagram map  is  
applied to a twisted polygon in 3D}
\label{T1-spiral}
\end{figure}

While this map is in a sense unique in 2D, its generalizations to higher dimensions seem to allow more freedom. A natural requirement for such generalizations, though,  is their integrability.
In  \cite{KS}  we observed that there is no natural generalization of this map to polyhedra and suggested a natural integrable generalization of the pentagram map to generic twisted space 
polygons (see Figure \ref{T1-spiral}). This generalization in any dimension was defined 
via intersections of ``short diagonal" hyperplanes, which are symmetric 
higher-dimensional analogs of polygon diagonals, see Section \ref{sect:any-diag} below.
This map turned out to be scale invariant (see
 \cite{OST99} for 2D, \cite{KS} for 3D, \cite{Beffa_scale} for higher D) and  integrable in any dimension as it admits a Lax representation with a spectral parameter  \cite{KS}.

 A  different integrable generalization to higher dimensions was proposed in \cite{GSTV}, where the pentagram map was defined not on generic, but on the so-called
  corrugated polygons. These are piecewise linear curves in $\RP^d$, whose pairs of edges with indices differing by $d$ lie in one and the same two-dimensional
   plane. It turned out that the pentagram map on corrugated polygons is integrable
 and it admits an explicit description of the Poisson structure, a cluster algebra structure, and other
 interesting features \cite{GSTV}.

 In this paper we present a  variety of integrable   generalized pentagram maps, which unifies these two approaches. ``Primary integrable maps" in our construction are called
 the  dented pentagram maps.
 These maps are defined for  generic twisted polygons in $\RP^d$.
 It turns out that the pentagram maps for corrugated polygons considered in \cite{GSTV} are a particular case (more precisely, a restriction) of these dented maps.
 We describe in detail how to perform such a reduction in Section \ref{S:corrug}.

To define the dented maps, we propose a definition of a ``dented diagonal hyperplane"
depending on a parameter $m=1,...,d-1$,  where $d$ is the dimension of the projective space.
The parameter $m$ marks the skipped vertex of the polygon,
and in dimension $d$ there are $d-1$ different dented integrable maps.
The vertices in the ``dented diagonal hyperplanes" are chosen in a non-symmetric way
(as opposed to the unique symmetric choice in \cite{KS}).
 We would like to stress that in spite of a non-symmetric choice, the integrability property is preserved,  and each of the dented maps can be regarded as a natural  generalization of the classical 2D pentagram map  of \cite{Schwartz}.
We describe the geometry and  Lax representations of the dented maps and their generalizations, the deep-dented pentagram maps, and prove their algebraic-geometric integrability in 3D.
In a sense, from now on a new  challenge might be 
to find examples of non-integrable Hamiltonian maps of pentagram type, cf. \cite{KS14}. 

We emphasize that often throughout the paper we understand {\it integrability} as
the existence of a Lax representation with a spectral parameter
corresponding to scaling invariance of a given dynamical system. 
We show how it is used to prove algebraic-geometric integrability for the primary maps in $\CP^3$.
In any dimension,  the Lax representation provides first integrals (as the coefficients 
of the corresponding spectral curve) and allows one to use algebraic-geometric 
machinery to prove various integrability properties. 
We also note that while most of the paper deals with $n$-gons satisfying the condition  
$gcd(n,d+1)=1$, the results hold in full generality and we show how they are adapted 
to the general setting in Section \ref{nonprimes}. 
While most of definitions below work both over $\R$ and $\C$, \ throughout the paper 
we describe the geometric features of pentagram maps over $\R$, while their 
Lax representations over $\C$.

\medskip

Here are the main results of the paper.
\smallskip

$\bullet$ We define generalized pentagram maps $T_{I,J}$ on (projective equivalence classes of)
twisted polygons in $\RP^d$, associated with $(d-1)$-tuple of numbers $I$ and $J$:
the tuple $I$ defines which vertices to take in the definition of the diagonal hyperplanes,
while the tuple $J$ determines which of the hyperplanes to intersect in order to get the image point.
In Section \ref{sect:any-diag} we prove the duality between such pentagram maps:
$$
T_{I,J}^{-1}=T_{J^*,I^*}\circ Sh\,,
$$
where $I^*$ and $J^*$ stand for the   $(d-1)$-tuples  taken in the opposite order and  $Sh$ is any shift in the indices  of polygon vertices.

\smallskip

$\bullet$  The {\it dented pentagram maps} $T_m$ on polygons  $(v_k)$ in $\RP^d$ are defined by intersecting $d$ consecutive diagonal hyperplanes. Each hyperplane $P_k$ passes through all vertices but one  from $v_k$ to $v_{k+d}$ by skipping only the vertex $v_{k+m}$.
The main theorem on  such maps  is the following (cf. Theorem \ref{thm:lax_anyD}):

\begin{theorem}
The dented  pentagram map $T_m$ on both twisted and closed 
$n$-gons in any dimension $d$ and any
$m=1,...,d-1$ is an integrable system in the sense that it admits a Lax representation
with a spectral parameter.
\end{theorem}

We also describe the dual dented maps, prove their scale invariance (see Section \ref{sect:dual}),
and study their geometry in detail. Theorem \ref{thm:comparison} shows that in
 dimension 3 the algebraic-geometric integrability follows from the proposed Lax representation for
 both dented pentagram maps and the short-diagonal pentagram map.

\medskip


$\bullet$   The continuous limit of any dented pentagram map $T_m$
(and more generally, of any generalized pentagram map) in dimension $d$  is the $(2, d+1)$-KdV flow of the Adler-Gelfand-Dickey  hierarchy on the circle, see Theorem \ref{thm:cont}.
For 2D this is the classical  Boussinesq equation  on the circle: $u_{tt}+2(u^2)_{xx}+u_{xxxx}=0$, which appears as the continuous limit of the 2D pentagram map \cite{OST99, Sch08}.

$\bullet$   Consider the space of  corrugated polygons in
$\RP^d$, i.e., twisted polygons, whose vertices $v_{k-1}, v_{k}, v_{k+d-1},$ and $v_{k+d}$
span a projective two-dimensional plane  for every $k\in \Z$, following \cite{GSTV}. It turns out that the pentagram
map $T_{cor}$ on them can be viewed as a particular case of the dented pentagram map, see
Theorem \ref{thm:restr_to_corr}:

\begin{theorem}
This pentagram map $T_{cor}$ is a restriction of the dented pentagram map $T_m$ for any $m=1,..., d-1$ from generic n-gons ${\mathcal P}_n$ in $\RP^d$ to corrugated ones ${\mathcal P}_n^{cor}$ (or differs from it by a shift in vertex indices). In particular, these restrictions for different $m$  coincide modulo an index shift.
\end{theorem}

We also describe the algebraic-geometric integrability for corrugated pentagram map in $\CP^3$, see Section \ref{corr3D}.
\medskip

$\bullet$   Finally, we provide an application of the use of dented pentagram maps.
The latter can be regarded as  ``primary" objects, simplest integrable systems of pentagram type.
By considering more general diagonal hyperplanes, such as ``deep-dented diagonals",
i.e., those skipping more than one vertex, one can construct new integrable systems,
see Theorem \ref{thm:ddd}:

\begin{theorem}
The deep-dented pentagram maps in $\RP^d$  are restrictions of integrable systems to  invariant submanifolds and have   Lax representations with a spectral parameter.
\end{theorem}

The main tool to prove integrability in this more general setting is an introduction of the corresponding notion
of {\it partially corrugated polygons},  occupying an intermediate position between corrugated
and generic ones, see Section \ref{sect:appl}. The pentagram map on such 
partially corrugated polygons also turns out to be integrable. 
This work brings about the following question, which 
manifests  the change of perspective on generalized pentagram maps:  

\begin{problem}
Is it possible to choose the diagonal hyperplane so that the corresponding pentagram map turned out to be {\rm non-integrable}?
\end{problem}

Some numerical evidence in this direction is presented in \cite{KS14}.

\bigskip

{\bf Acknowledgments}.
We are grateful to S.~Tabachnikov for useful discussions.
B.K. and F.S. were  partially supported by NSERC research grants.
B.K. is grateful to the Simons Center for Geometry and Physics in Stony Brook 
for support and hospitality;
F.S. acknowledges the support  of the Fields Institute in Toronto and the CRM in Montreal.

\bigskip


\section{Duality of pentagram maps in higher dimensions}\label{sect:any-diag}

We start with  the notion of a twisted $n$-gon in dimension $d$.

\begin{definition}\label{tw-ngon}
{\rm
A {\it twisted $n$-gon} in a projective space $\RP^d$ with a monodromy $M \in SL_{d+1}(\R)$
is a map $\phi: \Z \to \RP^d$, such that
$\phi(k+n) =  M \circ \phi(k)$ for each $k\in \Z$ and where $M$ acts naturally on $\RP^d$.
Two twisted $n$-gons are {\it equivalent} if there is a transformation $g \in SL_{d+1}(\R)$ such that $g \circ \phi_1=\phi_2$.
}
\end{definition}

We assume that the vertices $v_k:=\phi(k), \; k \in \Z,$ are in general position (i.e., no $d+1$ consecutive vertices lie in the same hyperplane in $\RP^d$), and denote by   ${\mathcal P}_n$ the space of generic twisted $n$-gons considered up to the above equivalence.
Define   general pentagram  maps as follows.

\begin{definition}\label{def:I-diag}
{\rm
Let $I=(i_1,...,i_{d-1})$ and $J=(j_1,...,j_{d-1})$ be two $(d-1)$-tuple of numbers $i_\ell, j_m\in \N$.
For a generic twisted $n$-gon in $\RP^d$ one can define a {\it $I$-diagonal hyperplane}  $P_k$
as the one passing through $d$ vertices of the $n$-gon by taking every $i_\ell$th vertex starting at the point $v_k$, i.e.,
$$
P_k:=(v_k, v_{k+i_1}, v_{k+i_1+i_2},..., v_{k+i_1+...+i_{d-1}})\,,
$$
see Figure \ref{fig:T312}.

\begin{figure}[hbtp!]
\centering
\includegraphics[width=2.9in]{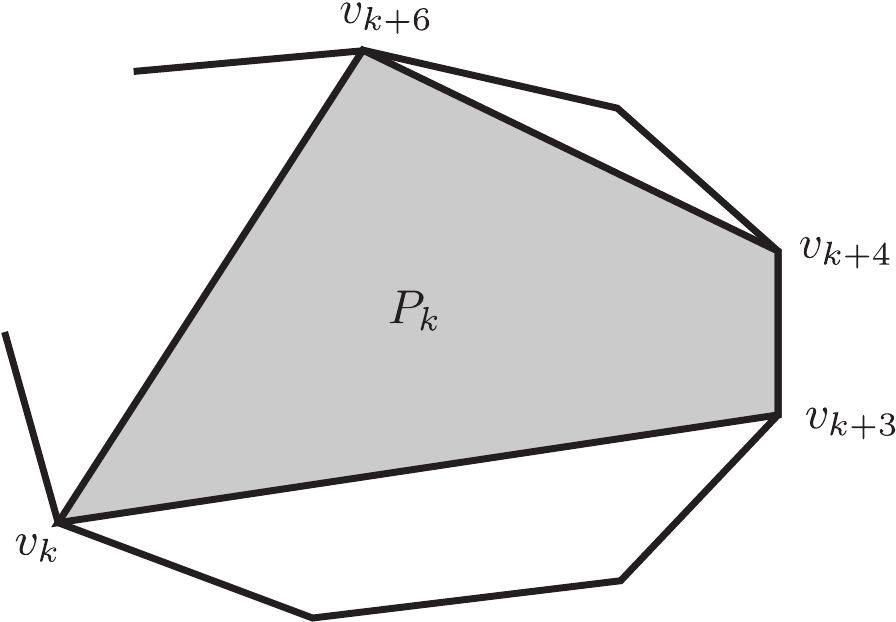}
\caption{\small The diagonal hyperplane for the jump tuple $I=(3,1,2)$ in $\RP^4$.} \label{fig:T312}
\end{figure}

The image of the vertex $v_k$ under the {\it
generalized pentagram map} $T_{I,J}$ is defined by intersecting
every $j_m$th out of the $I$-diagonal hyperplanes  starting with
$P_k$:
$$
T_{I,J}v_k:=P_{k}\cap P_{k+j_1}\cap P_{k+j_1+j_2}\cap...\cap P_{k+j_1+...+j_{d-1}}\,.
$$
(Thus $I$ defines the structure of the diagonal hyperplane, while $J$ governs which of them to intersect.)
The corresponding map $T_{I,J}$ is considered (and is generically defined) on the space ${\mathcal P}_n$  of equivalence classes of  $n$-gons in
$\RP^d$.
As usual, we assume that the vertices are in ``general position,'' and every $d$
 hyperplanes $P_i$ intersect at one point in $\RP^d$.
}
\end{definition}

\begin{example}
{\rm
Consider the case of $I=(2,2,...,2)$ and $J=(1,1,...,1)$ in $\RP^d$.
 This choice of $I$ corresponds to ``short diagonal hyperplanes",
  i.e., every $I$-diagonal hyperplane passes through $d$ vertices by taking every other vertex of the twisted polygon.
   The choice of $J$ corresponds to taking intersections of $d$ consecutive hyperplanes.
    This recovers the definition of the short-diagonal (or higher) pentagram  maps from \cite{KS}.
    Note that the classical 2D pentagram map
    has $I$ and $J$ each consisting of one number: $I=(2)$ and $J=(1)$.
}
\end{example}

Denote by $I^*=(i_{d-1},...,i_{1})$ the  $(d-1)$-tuple $I$ taken in the opposite order and
by $Sh$ the operation of any index shift on the sequence of vertices.

\begin{theorem}{\bf (Duality)}\label{thm:duality}
There is the following duality for the generalized pentagram maps $T_{I,J}$:
$$
T_{I,J}^{-1}=T_{J^*,I^*}\circ Sh\,,
$$
where $Sh$ stands for some shift in indices of vertices.
\end{theorem}

\proof
To prove this theorem we introduce the following duality maps, cf. \cite{OST99}.

\begin{definition}
{\rm
Given a generic sequence of points $\phi(j) \in \RP^d, \; j \in \Z,$  and
a $(d-1)$-tuple $I=(i_1,...,i_{d-1})$
we define the following {\it sequence of hyperplanes} in $\RP^d$:
$$
\alpha_I(\phi(j)):=(\phi(j), \phi(j+i_1),..., \phi(j+i_1+...+i_{d-1}))\,,
$$
which is regarded as a sequence of points in the dual space: $\alpha_I(\phi(j))\in (\RP^d)^*$.
}
\end{definition}

The generalized pentagram map $T_{I,J}$ can be defined as a composition of two such maps up to a shift of indices:
$T_{I,J}=\alpha_I\circ\alpha_J\circ Sh$.

Note that for a special $I=(p,p,...,p)$ the maps $\alpha_I$ are involutions modulo index shifts (i.e., $\alpha_I^2=Sh$), but for a general $I$ the  maps $\alpha_I$ are no longer involutions. However, one can see from their construction that they have the following duality property:
$\alpha_I\circ \alpha_{I^*}=Sh$ and they commute with index shifts: $\alpha_I\circ Sh=Sh\circ \alpha_I$.

Now we see that
$$
T_{I,J}\circ T_{J^*,I^*}=(\alpha_I\circ\alpha_J\circ Sh)\circ
(\alpha_{J^*}\circ\alpha_{I^*}\circ Sh)=Sh\,,
$$
as required.
\proofend

\begin{remark}\label{rem:T_pr}
{\rm
For $d$-tuples $I=(p,p,...,p)$ and $J=(r,r,...,r)$ the
generalized pentagram maps  correspond to the general pentagram maps $T_{p,r}=T_{I,J}$
discussed in \cite{KS}, and they possess the following duality: $T^{-1}_{p,r}=T_{r,p}\circ Sh$.

Note that in \cite{Beffa} one considered an intersection of the hyperplane $P_k$
with a chord joining two vertices, which leads to a different generalization of the pentagram map
and for which an analog of the above duality is yet unknown.

\smallskip

In the paper \cite{KS} we studied the case $T_{2,1}$ of short diagonal hyperplanes: $I=(2,2,...,2)$ and $J=(1,1,...,1)$, which is a very symmetric way of choosing the hyperplanes and their intersections.
In this paper we consider the general, non-symmetric choice of vertices.
}
\end{remark}

\begin{theorem}\label{thm:dual}
If $J=J^*$  (i.e., $\alpha_J$ is an involution), then modulo a shift in indices

$i)$ the pentagram maps $T_{I,J}$
and $T_{J, I^*}$ are inverses to each other;

$ii)$ the pentagram maps $T_{I,J}$ and $T_{J, I}$ (and hence $T_{I,J}$
and $T_{I^*, J}^{-1}$) are conjugated  to each other, i.e.,
the map  $\alpha_J$ takes the map $T_{I,J}$ on $n$-gons in $\RP^d$ into the map
 $T_{J,I}$ on $n$-gons in $(\RP^d)^*$.

 In particular, all four maps $T_{I,J}, T_{I^*, J},
 T_{J, I}$ and $T_{J, I^*}$ are integrable or non-integrable simultaneously. Whenever they are integrable,  their integrability characteristics, e.g., the dimensions of invariant tori, periods of the corresponding orbits, etc., coincide.
\end{theorem}

\proof
The statement $i)$ follows from Theorem \ref{thm:duality}. To prove $ii)$ we note that
for $J=J^*$ one has $\alpha_J^2=Sh$ and therefore
$$
\alpha_J\circ T_{I, J} \circ \alpha_J^{-1}=\alpha_J\circ (\alpha_I\circ\alpha_J\circ Sh)\circ \alpha_J
=(\alpha_J\circ \alpha_I\circ Sh)\circ\alpha_J^2=T_{J,I}\circ Sh\,.
$$
Hence  modulo index shifts, the pentagram map $T_{I, J}$ is conjugated to $T_{J,I}$, while by the statement $i)$ they are also inverses of  $T_{J, I^*}$ and $T_{I^*, J}$ respectively. This proves the theorem.
\proofend


\section{Dented pentagram maps}

\subsection{Integrability of dented pentagram maps}\label{sect:dent}

From now on we consider the case of $J={\mathbf 1}:=(1,1,...,1)=J^*$ for different $I$'s, i.e., we take the intersection of {\it consecutive} $I$-diagonal hyperplanes.

\begin{figure}[hbtp!]
\centering
\includegraphics[width=2.9in]{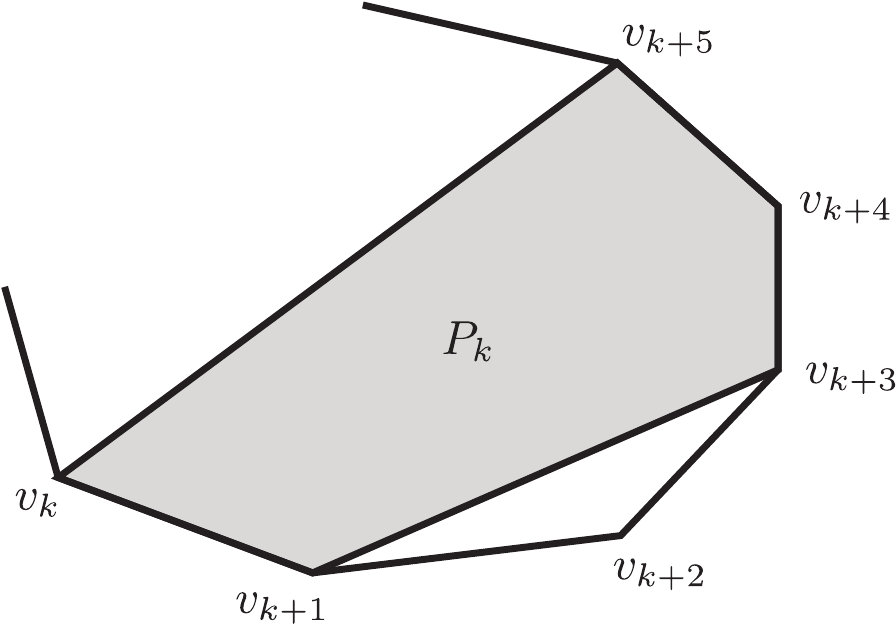}
\caption{\small The dented diagonal hyperplane $P_k$ for $m=2$ in $\RP^5$.}
\label{fig:dented-plane}
\end{figure}

\begin{definition}
{\rm
Fix an integer parameter $m\in \{1,...,d-1\}$ and for the $(d-1)$-tuple $I$
we set  $I=I_m:=(1,...,1,2,1,...,1)$, where
the only value 2 is situated at the $m$th place: $i_m=2$ and $i_\ell=1$ for $\ell\not=m$.
This choice of the tuple $I$ corresponds to the  diagonal plane $P_k$ which passes through
consecutive vertices $v_k, v_{k+1},...,v_{k+m-1}$, then skips vertex $v_{k+m}$
and continues passing through consecutive vertices $v_{k+m+1},...,v_{k+d}$:
$$
P_k:=(v_k, v_{k+1},...,v_{k+m-1},v_{k+m+1},v_{k+m+2},...,v_{k+d})\,.
$$
We call such a plane $P_k$ a {\it dented} (or {\it $m$-dented}) {\it diagonal plane},
as it is ``dented" at the vertex $v_{k+m}$, see Figure \ref{fig:dented-plane}.
We define the {\it dented pentagram map} $T_m$ by  intersecting $d$ consecutive planes
$P_k$:
$$
T_m v_k:=P_{k}\cap P_{k+1}\cap ...\cap P_{k+d-1}\,.
$$
In other words,  the dented pentagram map is $T_m:=T_{I_m,{\mathbf 1}}$, i.e. $ T_{I_m,J}$ where
$J={\mathbf 1}$.
}
\end{definition}

\begin{corollary}\label{cor:T_m}
The dented map  $T_m$ is conjugated (by the involution $\alpha_{\mathbf 1}$)
to $T^{-1}_{d-m}$  modulo shifts.
\end{corollary}

\proof
Indeed,  $I_m=I^*_{d-m}$ and hence, due to Theorem \ref{thm:dual}, one has
$\alpha_{\mathbf 1}\circ T_m\circ \alpha_{\mathbf 1} =\alpha_{\mathbf 1}\circ T_{I_m,{\mathbf 1}}\circ \alpha_{\mathbf 1}
=T_{{\mathbf 1},I_m}\circ Sh
=T^{-1}_{I^*_m,{\mathbf 1}}\circ Sh=T^{-1}_{I_{d-m},{\mathbf 1}}\circ Sh=T^{-1}_{d-m}\circ Sh$, where $\alpha_{\mathbf 1}$ stands for $\alpha_J$ for $J=(1,...,1)$.
\proofend

One can also see that for $m=0$ or $m=d$ all the vertices defining the hyperplane $P_k$ are taken consecutively, and the corresponding map $T_m$ is the identity modulo a shift in indices of $v_k$.

\medskip

For 2D the only option for the dented map is $I=(2)$ and $J=(1)$, and where $m=1$,
so the corresponding map $T_m$ coincides with the classical
pentagram transformation $T=T_{2,1}$ in 2D.
Thus the above definition of maps $T_m$ for various $m$  is another natural
higher-dimensional generalization of the 2D pentagram map.
Unlike the definition of the short-diagonal pentagram map $T_{2,1}$ in $\RP^d$, the dented pentagram map
 is not unique for each dimension $d$, but also has one more integer parameter $m=1,...,d-1$.

\bigskip

It turns out that the dented pentagram map $T_m$ defined this way, i.e., defined as $T_{I_m,{\mathbf 1}}$
for $I_m=(1,...,1,2,1,...,1)$ and
${\mathbf 1}=(1,1,...,1)$, has a special scaling invariance. To describe it we need
to introduce coordinates on the space ${\mathcal P}_n$ of twisted $n$-gons.

Now we complexify the setting and consider the spaces and
maps over $\C$.

\begin{remark-definition}\label{diff-eq}
{\rm
One can show that there exists a  lift of the vertices $v_k=\phi(k) \in \CP^d$
to the vectors $V_k \in \C^{d+1}$ satisfying $\det(V_j, V_{j+1}, ..., V_{j+d})=1$ and $ V_{j+n}=MV_j,\; j \in \Z,$ where
$M\in SL_{d+1}(\C)$, provided that the condition $gcd(n,d+1)=1$ holds.
The corresponding lifted vectors satisfy the difference equations have the form
\begin{equation}\label{eq:difference_anyD}
V_{j+d+1} = a_{j,d} V_{j+d} + a_{j,d-1} V_{j+d-1} +...+ a_{j,1} V_{j+1} +(-1)^{d} V_j,\quad j \in \Z,
\end{equation}
with $n$-periodic coefficients in the index $j$.
This allows one to introduce {\it coordinates} $\{ a_{j,k} ,\;0\le j\le n-1, \; 1\le k\le d \}$ on the space of twisted $n$-gons in $\CP^d$. In the theorems below we assume the condition $gcd(n,d+1)=1$
whenever we use explicit formulas in the coordinates $\{ a_{j,k}\}$. However
the statements  hold in full generality and we discuss how the corresponding formulas are being adapted in Section \ref{nonprimes}. 
(Strictly speaking, the  lift from vertices to vectors is not unique, because it is defined up to simultaneous 
multiplication of all vectors by $\varepsilon$, where $\varepsilon^{d+1}=1$, 
but  the coordinates $\{ a_{j,k}\}$ are well-defined as they have the same values for all lifts.)\footnote{Note also 
that over $\R$  for odd $d$  to obtain the lifts of $n$-gons from $\RP^d$ to $\R^{d+1}$ 
one might need to switch the sign of the monodromy matrix: $M \to -M \in SL_{d+1}(\R)$, since the field is not algebraically closed. These monodromies in $SL_{d+1}(\R)$ correspond to the same projective monodromy
in $PSL_{d+1}(\R)$.}
}
\end{remark-definition}

\begin{theorem} {\bf (Scaling invariance)}\label{thm:scaling}
The dented pentagram map $T_m$ on twisted $n$-gons in $\CP^d$ with hyperplanes $P_k$ defined by taking the vertices in a row but skipping the $m${\rm th} vertex is invariant with respect to the following scaling transformations:
$$
a_{j,1} \to s^{-1}a_{j,1} ,\; a_{j,2} \to s^{-2}a_{j,2} ,\;...\:, \; a_{j,m} \to s^{-m} a_{j,m},
$$
$$
 a_{j,m+1} \to s^{d-m}a_{j,m+1}, \;...\:, \;
a_{j,d} \to s a_{j,d}
$$
for all $s\in \C^*$.
\end{theorem}

For $d=2$ this is the case of the classical pentagram map, see \cite{OST99}.
We prove this theorem in Section \ref{sect:scale_proof}.
The above scale invariance implies the Lax representation, which
opens up the possibility to establish algebraic-geometric
integrability of the dented pentagram maps.

 \begin{remark}
 {\rm
Recall that a discrete Lax equation with a spectral
 parameter is a representation of a dynamical system in the form
\begin{equation}\label{lax-eq}
 L_{j,t+1}(\lambda) = P_{j+1,t}(\lambda) L_{j,t}(\lambda) P_{j,t}^{-1}(\lambda),
 \end{equation}
where $t$ stands for the discrete time variable,  $j$ refers to the vertex index,
and $\lambda$ is a complex spectral parameter. It is a discrete version
of the classical zero curvature equation $\partial_tL-\partial_xP=[P,L]$.
}
 \end{remark}

\begin{theorem}{\bf (Lax form)}\label{thm:lax_anyD}
The dented  pentagram map $T_m$ on both twisted and closed 
$n$-gons in any dimension $d$ and any
$m=1,...,d-1$ is an integrable system in the sense that it admits a Lax representation
 with a spectral parameter. In particular, for $gcd(n,d+1)=1$ the Lax matrix is 
\[
L_{j,t}(\lambda) =
\left(
\begin{array}{cccc|c}
0 & 0 & \cdots & 0    &(-1)^d\\ \cline{1-5}
\multicolumn{4}{c|}{\multirow{4}*{$D(\lambda)$}} & a_{j,1}\\
&&&& a_{j,2}\\
&&&& \cdots\\
&&&& a_{j,d}\\
\end{array}
\right)^{-1},
\] 
with the diagonal $(d \times d)$-matrix $D(\lambda)={\rm diag}(1,...,1,\lambda, 1,...1)$, where
the spectral parameter $\lambda$ is situated at the $(m+1)${\rm th}  place, and an appropriate matrix $P_{j,t}(\lambda)$.
\end{theorem}


\proof
Rewrite the difference equation  \eqref{eq:difference_anyD} in the matrix form.
It is equivalent to the relation $(V_{j+1},V_{j+2},...,V_{j+d+1})=(V_j,V_{j+1},...,V_{j+d})N_{j}$, where
the transformation matrix $N_{j}$ is
\[
N_{j} :=
\left(
\begin{array}{ccc|c}
0 &  \cdots & 0    &(-1)^d\\ \cline{1-4}
\multicolumn{3}{c|}{\multirow{3}*{\rm{Id}}} & a_{j,1}\\
&&& \cdots\\
&&& a_{j,d}\\
\end{array}
\right),
\] 
and  where $\mathrm{Id}$ stands for the identity $(d\times d)$-matrix.

It turns out that the monodromy  $M$ for  twisted $n$-gons is always conjugated to
the product $\tilde{M}:=N_{0} N_{1}...N_{n-1}$, see Remark \ref{rem:monodromy} below.
Note that the pentagram map defined
on classes of projective equivalence preserves the conjugacy class of $M$ and hence that of
 $\tilde{M}$. Using the scaling invariance of the pentagram map $T_m$, replace
$a_{j,k}$ by $s^* a_{j,k}$ for all $k$ in the right column of $N_j$ to obtain a new matrix $N_j(s)$.
 The  pentagram map preserves the conjugacy class of  the new monodromy $\tilde{M}(s):=N_0(s) ...N_{n-1}(s)$ for any $s$,
 that is, the monodromy can only change to a conjugate one during its pentagram evolution:
$ \tilde{M}_{t+1}(s) = P_{t}(s) \tilde{M}_t(s) P_{t}^{-1}(s)$. Then $N_j(s)$ (or, more precisely, $N_{j,t}(s)$ to emphasize its dependence on $t$), being a discretization of the monodromy
$\tilde{M}$, could be
taken as a Lax matrix $L_{j,t}(s)$. The gauge transformation 
$L_{j,t}^{-1}(\lambda) := \left(g^{-1} N_j(s) g \right)/s$ for
$g = \text{diag}(s^{-1},s^{-2},...,s^{-m-1},s^{d-m-1},...,s,1)$ and $\lambda \equiv s^{-d-1}$ simplifies the formulas  and gives the required matrix $L_{j,t}(\lambda)$.

Closed polygons are subvarieties defined by polynomial relations on 
coefficients $a_{j,k}$. These relations ensure that the monodromy $\tilde{M}(s)$
has an eigenvalue of multiplicity $d+1$ at $s=1$. 
\proofend

\begin{remark}\label{rem:monodromy}
{\rm Define the current monodromy  $\tilde M_j$ for  twisted $n$-gons by the relation
$$
(V_{j+n},V_{j+n+1},...,V_{j+n+d})=(V_j,V_{j+1},...,V_{j+d})\tilde M_{j},
$$
i.e., as the product $\tilde{M}_j:=N_{j} N_{j+1}...N_{j+n-1}$.
Note that $\tilde{M}_j$ acts on matrices by multiplication on the right,
whereas in Definition~\ref{tw-ngon} the monodromy $M$ acts on vectors $V_j$ on the left.
The theorem above uses the following fact:
}
\end{remark}

\begin{lemma}
All current monodromies $\tilde M_j$ lie in the same conjugacy class in $SL_{d+1}(\C)$ as $M$.
\end{lemma}

\proof
All products $\tilde{M}_j:=N_{j} N_{j+1}...N_{j+n-1}$ are conjugated:
$\tilde{M}_{j+1}=N_{j}^{-1}\tilde{M}_j N_j$ for all $j\in \Z$, since $N_j=N_{j+n}$.
Furthermore,
$$
(V_j,V_{j+1},...,V_{j+d})\tilde M_{j}(V_j,V_{j+1},...,V_{j+d})^{-1}\-=(V_{j+n},V_{j+n+1},...,V_{j+n+d})(V_j,V_{j+1},...,V_{j+d})^{-1}
$$
$$
= M(V_j,V_{j+1},...,V_{j+d})(V_j,V_{j+1},...,V_{j+d})^{-1}=M\,.
$$
\proofend

To prove the scale invariance of dented pentagram maps we need to introduce the notion of  the corresponding dual  map.

 \begin{definition}
 {\rm
The {\it dual dented pentagram map} $\widehat{T}_m$ for twisted polygons in $\RP^d$ or $\CP^d$
is defined as $\widehat{T}_m:=T_{{\mathbf 1},I^*_m}$ for  $I^*_m=(1,...,1,2,1,...,1)$ where
2 is at the $(d-m)$th place and ${\mathbf 1}=(1,1,...,1)$. In this case the diagonal planes $P_k$ are
defined by taking $d$ consecutive vertices of the polygon starting with the vertex $v_k$,
but to define the image $\widehat T_mv_k$  of the vertex $v_k$ one takes the intersection
$P_k\cap P_{k+1}\cap ...\cap P_{k+d-m-1}\cap P_{k+d-m+1}\cap ...\cap P_{k+d}$ of all but one consecutive planes by skipping only the plane $P_{k+d-m}$.
}
\end{definition}

\smallskip

 \begin{remark}
{\rm
According to Theorem \ref{thm:dual},
the dual map satisfies $\widehat{T}_m={T}^{-1}_m\circ Sh$.
In particular, the dual map $\widehat{T}_m$ is also integrable
and has the same scaling properties and the Lax matrix as ${T}_m$.
The dynamics for $\widehat{T}_m$ is obtained by reversing time in the dynamics of ${T}_m$.
Moreover, the map $\widehat{T}_m$ is conjugated to $T_{d-m}$ (modulo shifts) by means of the involution $\alpha_{\mathbf 1}$.
}
 \end{remark}


\begin{example}
{\rm
In dimension $d=3$ one has the following explicit Lax representations.
For the case of $T_1$ (i.e., $m=1$) one sets $D(\lambda)=(1,\lambda,1)$.
 The dual map $\widehat{T}_1$, being the inverse of $T_1$, has the same Lax  form and scaling.

For the map $T_2$ (where $m=2$) one has $D(\lambda)=(1,1,\lambda)  $. Similarly,
 $\widehat{T}_2$ is the inverse of $T_2$. Note that the maps $T_1$ and $T_2^{-1}$ are conjugated to each other by means of the involution $\alpha_{\mathbf 1}$ for ${\mathbf 1}=(1,1)$. They have the same  dimensions of invariant tori, but their Lax forms differ.
 }
\end{example}

\begin{example}
{\rm
In dimension  $d=4$ one has two essentially different cases, according to whether the dent is on a side of the diagonal plane or in its middle. Namely,  the map $T_2$ is the case where the diagonal hyperplane is dented in the middle point, i.e., $m=2$ and $I_m=(1,2,1)$.
In this case $D(\lambda)=(1,1,\lambda,1)$.

For the side case consider the map $T_1$ (i.e., $m=1$ and $I_m=(2,1,1)$), where $D(\lambda)=(1,\lambda,1,1)$.
The dual map $\widehat{T}_1$ is the inverse of $T_1$ and has the same Lax form.
The map $T_3$ has the  Lax form with $D(\lambda)=(1,1,1,\lambda)$   and
 is conjugate to the inverse  $T_1^{-1}$, see Corollary \ref{cor:T_m}.
}
\end{example}


\subsection{Coordinates in the general case}\label{nonprimes}

In this section we describe how to introduce coordinates on the space of twisted polygons
for any $n$. If  $gcd(n,d+1)\not=1$ one can use quasiperiodic coordinates $a_{j,k}$
subject to a certain equivalence relation, instead of periodic ones, cf. Section 5.3 in \cite{KS}.

\begin{definition}\label{def:quasi-abc}
{\rm
Call $d$ sequences of coordinates $\{a_{j,k}, k=1,...,d, \;j\in \Z\}$ $n$-{\it quasi-periodic} if
there is a $(d+1)$-periodic sequence $t_j,\; j \in \Z$, satisfying
$t_j t_{j+1} ... t_{j+d}=1$  and such that $a_{j+n, k} = a_{j,k}\cdot{t_j}/{t_{j+k}}$
for each $j \in \Z$.
}
\end{definition}

This definition arises from the fact that there are different lifts of vertices 
$v_j\in \CP^d$ to vectors $V_j\in \C^{d+1}, \; j \in \Z,$ so that 
$\det(V_j, V_{j+1}, ..., V_{j+d})=1$ and $ v_{j+n}=Mv_j$ for $M\in SL_{d+1}(\C)$ and $j \in \Z$. (The later monodromy condition on vertices $v_j$ is weaker than the condition $ V_{j+n}=MV_j$ on lifted vectors  in Definition \ref{diff-eq}.)
We take arbitrary lifts $V_0,  ..., V_{d-1}$
of the first $d$ vertices $v_0, ..., v_{d-1}$  and then obtain that 
$ V_{j+n}=t_jMV_j$, where $t_j t_{j+1}\dots, t_{j+d}=1$ and $t_{j+d+1}=t_j $ for all $ j \in \Z$,
see details in \cite{OST99, KS}. 
This way twisted $n$-gons are described by quasiperiodic coordinate sequences $a_{j,k}, \; k=1,..., d, \; j\in \Z$, with the equivalence furnished by different choices of $t_j,\; j \in \Z$.  Indeed, 
the defining relation (\ref{eq:difference_anyD}) after adding $n$ to all indices $j$'s becomes the relation 
$$
t_jV_{j+d+1} = a_{j,d} V_{j+d} t_{j+d}
+ a_{j,d-1} V_{j+d-1} t_{j+d-1}+...+ a_{j,1} V_{j+1} t_{j+1}+(-1)^{d} V_j t_j,\quad j \in \Z,
$$
which is consistent with the quasi-periodicity condition on $\{a_{j,k}\}$.

In the case when $n$ satisfies $gcd(n,d+1)=1$,  one can choose the parameters $t_j$ in such a way that the sequences  $\{a_{j,k}\}$ are $n$-periodic in $j$. 
For a general $n$, from $n$-quasi-periodic sequences  $\{a_{j,k}, k=1,...,d,\; j\in \Z\}$
one can construct $n$-periodic ones  (in $j$) as follows:
$$
\tilde a_{j,k}=\dfrac{a_{j+1,k-1}}{a_{j,k}a_{j+1,d}}
$$
for $j\in \Z$ and $k=1,...,d$, where one sets $a_{j,0}=1$ for all $j$. These new $n$-periodic coordinates   $\{\tilde a_{j,k} ,\;0\le j\le n-1, \; 1\le k\le d \}$
 are well-defined coordinates on twisted $n$-gons in $\CP^d$ (i.e., they do not depend on the choice of lift coefficients $t_j$). The periodic coordinates $\{\tilde a_{j,k}  \}$ are analogs of the cross-ratio coordinates $x_j, y_j$ in \cite{OST99} and $x_j, y_j, z_j$ in \cite{KS}.

\begin{theorem}{\bf (=\ref{thm:lax_anyD}$'$)}\label{thm:nonprimes}
The dented  pentagram map $T_m$ on $n$-gons in any dimension $d$ and any
$m=1,...,d-1$ is an integrable system. In the coordinates   $\{\tilde a_{j,k} \}$ its  Lax matrix is 
\[
\widetilde L_{j,t}(\lambda) =
\left(
\begin{array}{cccc|c}
0 & 0 & \cdots & 0    &(-1)^d\\ \cline{1-5}
\multicolumn{4}{c|}{\multirow{4}*{$A(\lambda)$}} & 1\\
&&&& 1\\
&&&& \cdots\\
&&&& 1\\
\end{array}
\right)^{-1},
\] 
where $A(\lambda)={\rm diag}(\tilde a_{j,1}, ..., \tilde a_{j,m},\lambda\tilde a_{j,m+1},\tilde a_{j,m+2},..., \tilde a_{j,d})$.
\end{theorem}

Note that  Lax matrices $\widetilde L$ and $L$ 
are related as follows: $\widetilde L_{j,t}(\lambda)=
a_{j+1, d}(h^{-1}_{j+1}  L_{j,t}(\lambda) h_j)$ for the matrix $h_j=\diag (1, a_{j,1},  a_{j,2},...,  a_{j,d})$.


\subsection{Algebraic-geometric integrability of pentagram maps in 3D}\label{sect:ag-in}

\begin{figure}[hbtp!]
\centering
\includegraphics[width=7in]{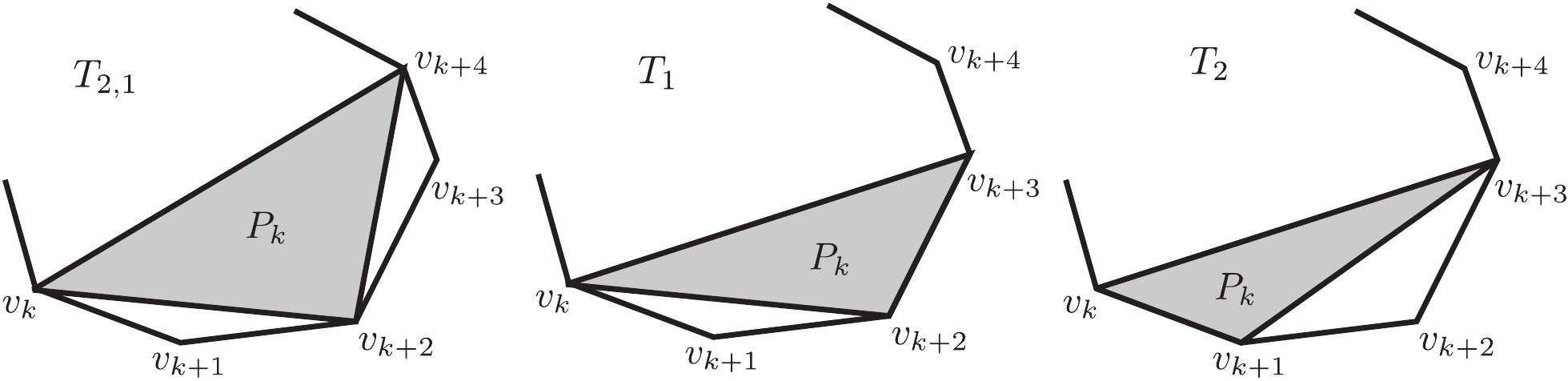}
\caption{\small Different diagonal planes in 3D: for $T_{2,1}, T_1,$ and $ T_2$.}
\label{fig:T-3D}
\end{figure}

The key ingredient responsible for algebraic-geometric integrability of the pentagram maps
is a Lax representation with a spectral parameter.
It allows one to construct the direct and inverse spectral transforms,
which imply that the dynamics of the maps
takes place on invariant tori, the Jacobians of the corresponding spectral curves.
The proofs in 3D case for the short-diagonal pentagram map $T_{2,1}$   are
presented in detail in \cite{KS} (see also \cite{FS} for the 2D case).
In dimension 3 we consider two dented pentagram maps $T_1$ and $T_2$,
where the diagonal hyperplane $P_k$ is dented on {\it different sides}
as opposed to the short-diagonal pentagram map $T_{2,1}$, where the diagonal hyperplane is dented on {\it both sides}, see Figure \ref{fig:T-3D}.

The proofs for the maps $T_1$ and $T_2$ follow the same line as in \cite{KS},
so in this section we present only the main statements and
outline the necessary changes.

 For simplicity, in this section we assume that $n$ is odd, which is equivalent 
 to the condition $gcd(n,d+1)=1$ for $d=3$
 (this condition may not appear for a different choice of coordinates, but the results of \cite{KS} show that the dimensions of tori  may depend on the parity of $n$). In this section 
 we  consider twisted polygons in the complex space $\CP^3$.

\begin{theorem}\label{thm:comparison}
In dimension 3  the dented pentagram maps on twisted $n$-gons
generically are fibered into (Zariski open subsets of) tori of  dimension 
$3\lfloor n/2\rfloor-1$ for $n$ odd and divisible by 3 and of  dimension  
$3\lfloor n/2\rfloor$ for $n$ odd and not divisible by 3.
\end{theorem}

Recall that for the short-diagonal pentagram map in 3D  the torus dimension is equal to 
$3\lfloor n/2\rfloor$ for any odd $n$, see \cite{KS}.

\proof
To prove this theorem we need the notion of a spectral curve. Recall that the product of Lax functions $L_j(\lambda),\; 0 \le j \le n-1,$ gives the monodromy operator $T_0(\lambda)$,
which determines the spectral function $R(k,\lambda):=\det{(T_0(\lambda) - k \,\text{Id})}$.
The zero set of $R(k,\lambda)=0$
is an algebraic curve in $\C^2$. A standard procedure (of adding the infinite points and normalization with
a few blow-ups) makes it into a compact Riemann surface, which we call the {\it spectral curve} and denote by $\Gamma$.
Its genus equals the dimension of the corresponding complex torus, its Jacobian and
Proposition \ref{prop:Jacobians} below shows how to find this genus.

As it is always the case with integrable systems, the spectral curve $\Gamma$
is an invariant of the map and the dynamics takes place on its Jacobian.
To describe the dynamics one introduces a {\it Floquet-Bloch solution} which
is formed by eigenvectors of the monodromy operator $T_0(\lambda)$.
After a certain normalization it becomes a uniquely defined meromorphic
vector function $\psi_0$ on  the spectral curve $\Gamma$.
Other Floquet-Bloch solutions are defined as the vector functions
$\psi_{i+1}=L_i...L_1 L_0 \psi_0,\; 0 \le i \le n-1$.
Theorem \ref{thm:comparison} is based on the study of $\Gamma$ and Floquet-Bloch solutions,
which we summarize in the tables below.

\smallskip

In each case, the analysis starts with an evaluation of the spectral function $R(k,\lambda)$.
Then we provide Puiseux series for the singular points at $\lambda=0$ and at $\lambda=\infty$.
They allow us to find the genus of the spectral curve and
the symplectic leaves for the corresponding Krichever-Phong's universal formula.
Then we describe the divisors of the Floquet-Bloch solutions, which are essential for  
constructing  the inverse spectral transform.

We start with reproducing the corresponding results for the short-diagonal 
map $T_{2,1}$ for odd $n$,
obtained in \cite{KS}.
We set $q:=\lfloor n/2 \rfloor.$ The tables below contain the information on the
Puiseux series of the spectral curve, Casimir functions of the pentagram dynamics, and divisors 
$(\psi_{i,k})$ of the components of
the Floquet-Bloch solutions $\psi_i$ (we refer to \cite{KS} for more detail).

\smallskip

\bull
{\rm
For the (symmetric) pentagram map $T_{2,1}$ defined in \cite{KS} by means of short diagonal hyperplanes, we have $D(\lambda)=\diag(\lambda,1,\lambda);$
$$
R(k,\lambda) = k^4 - \dfrac{k^3}{\lambda^n} \left( \sum_{j=0}^q G_j \lambda^j \right) + \dfrac{k^2}{\lambda^{q+n}} \left( \sum_{j=0}^q J_j \lambda^j \right)
- \dfrac{k}{\lambda^{2n}} \left( \sum_{j=0}^q I_j \lambda^j \right) + \dfrac{1}{\lambda^{2n}} = 0;
$$
\begin{center}
\begin{tabular}{|c|c|}
\hline
$\lambda=0$                                             & $\lambda=\infty$\\ \hline
$O_1: k_1 = 1/I_0 + {\mathcal O}(\lambda)$              &$W_{1,2}: k_{1,2,3,4}=k_\infty \lambda^{-n/2}(1+ {\mathcal O}(\lambda^{-1})),$\\
$O_2: k_{2,3}=\pm \sqrt{-I_0/G_0}\lambda^{-n/2} (1+{\mathcal O}(\sqrt{\lambda}))$ & where $k_\infty^4 + J_q k_\infty^2 + 1 = 0$.\\
$O_3: k_4 = G_0 \lambda^{-n}(1 + {\mathcal O}(\lambda))$ &\\ \hline
\multicolumn{2}{|c|}{$g = 3q$, the Casimirs are $I_0:=\prod_{j=0}^{n-1} a_{j,3};J_q;G_0:=\prod_{j=0}^{n-1} a_{j,1}$.}\\ \hline
\multicolumn{2}{|c|}{$(\psi_{i,1}) \ge -D +O_2 -i (O_2+O_3)+(i+1)(W_1+W_2)$}\\
\multicolumn{2}{|c|}{$(\psi_{i,2}) \ge -D +(1-i)(O_2+O_3)+i(W_1+W_2)$}\\
\multicolumn{2}{|c|}{$(\psi_{i,3}) \ge -D -i (O_2+O_3)+(i+1)(W_1+W_2)$}\\
\multicolumn{2}{|c|}{$(\psi_{i,4}) \ge -D +O_2+(1-i) (O_2+O_3)+i(W_1+W_2)$}\\ \hline

\end{tabular}
\end{center}
}

\bull
{\rm
 For the dented map $T_1$, we have $D(\lambda)=\diag(1,\lambda,1);$
$$
R(k,\lambda) = k^4 - \dfrac{k^3}{\lambda^q} \left( \sum_{j=0}^q G_j \lambda^j \right) + \dfrac{k^2}{\lambda^n} \left( \sum_{j=0}^{\lfloor 2n/3 \rfloor} J_j \lambda^j \right)
- \dfrac{k}{\lambda^n} \left( \sum_{j=0}^{\lfloor n/3 \rfloor} I_j \lambda^j \right) + \dfrac{1}{\lambda^n} = 0;
$$
\begin{center}
\begin{tabular}{|c|c|}
\hline
\multicolumn{2}{|c|}{\text{$n$ odd and not divisible by 3:} $n=6l+1$\text{ or } $n=6l+5$}\\ \hline
$\lambda=0$                                                                   & $\lambda=\infty$\\ \hline
$O_{1,2}: k_{1,2} = c_0 + {\mathcal O}(\lambda),$                             & $W_1: k_1=G_q+{\mathcal O}(\lambda^{-1}),$\\
where $c_0^2 J_0-c_0 I_0+1=0.$                                                & $W_2: k_{2,3,4}=G_q^{-1/3} \lambda^{-n/3}(1+ {\mathcal O}(\lambda^{-1/3}))$.\\
$O_3: k_{3,4}=\pm \sqrt{-J_0}\lambda^{-n/2} (1+{\mathcal O}(\sqrt{\lambda}))$ &\\ \hline
\multicolumn{2}{|c|}{$g = 3q$, the Casimirs are $I_0;J_0:=(-1)^n\prod_{j=0}^{n-1} a_{j,2};G_q:=\prod_{j=0}^{n-1} a_{j,1}$.}\\ \hline
\multicolumn{2}{|c|}{$(\psi_{i,1}) \ge -D + O_3 + 2W_2 + i (W_2-O_3)$}\\
\multicolumn{2}{|c|}{$(\psi_{i,2}) \ge -D + W_1 + 2W_2 + i (W_2-O_3)$}\\
\multicolumn{2}{|c|}{$(\psi_{i,3}) \ge -D + 2O_3 + i (W_2-O_3)$}\\
\multicolumn{2}{|c|}{$(\psi_{i,4}) \ge -D + 2O_3 + W_2 + i (W_2-O_3)$}\\ \hline
\end{tabular}
\end{center}

\begin{center}
\begin{tabular}{|c|c|}
\hline
\multicolumn{2}{|c|}{\text{$n$ odd and  divisible by 3:} $n=6l+3$}\\ \hline
$\lambda=0$                                                                   & $\lambda=\infty$\\ \hline
$O_{1,2}: k_{1,2} = c_0 + {\mathcal O}(\lambda),$                             & $W_1: k_1=G_q+{\mathcal O}(\lambda^{-1}),$\\
where $c_0^2 J_0-c_0 I_0+1=0.$                 & 
  $W_{2,3,4}: k_{2,3,4}=k_\infty \lambda^{-n/3}+ {\mathcal O}(\lambda^{-1})$, where\\
$O_3: k_{3,4}=\pm \sqrt{-J_0}\lambda^{-n/2} (1+{\mathcal O}(\sqrt{\lambda}))$ &
 $G_q k_\infty^3-J_{\lfloor 2n/3\rfloor} k_\infty^2 +  I_{\lfloor n/3\rfloor} k_\infty - 1 =0.$
\\ \hline
\multicolumn{2}{|c|}{$g = 3q-1$, the Casimirs are $I_0;J_0;G_q; J_{\lfloor 2n/3\rfloor}; 
I_{\lfloor n/3\rfloor}$.}\\ \hline
\end{tabular}
\end{center}
}

\bull
{\rm
For the dented map $T_2$, we have $D(\lambda)=\diag(1,1,\lambda);$
$$
R(k,\lambda) = k^4 - \dfrac{k^3}{\lambda^{\lfloor n/3 \rfloor}} \left( \sum_{j=0}^{\lfloor n/3 \rfloor} G_j \lambda^j \right) + \dfrac{k^2}{\lambda^{\lfloor 2n/3 \rfloor}} \left( \sum_{j=0}^{\lfloor 2n/3 \rfloor} J_j \lambda^j \right)
- \dfrac{k}{\lambda^n} \left( \sum_{j=0}^q I_j \lambda^j \right) + \dfrac{1}{\lambda^n} = 0\,.
$$
The analysis of the spectral curve proceeds similarly:
\begin{center}
\begin{tabular}{|c|c|}
\hline
\multicolumn{2}{|c|}{\text{$n$ odd and not divisible by 3:} $n=6l+1$\text{ or } $n=6l+5$}\\ \hline
$\lambda=0$                                                           & $\lambda=\infty$\\ \hline
$O_1: k_1 = 1/I_0 + {\mathcal O}(\lambda)$                                 & $W_{1,2}: k_{1,2}=c_1+{\mathcal O}(\lambda^{-1})$,\\
$O_2: k_{2,3,4}=I_0^{1/3} \lambda^{-n/3} (1+ {\mathcal O}(\lambda^{1/3}))$ & where $c_1^2-c_1 G_{\lfloor n/3 \rfloor}+J_{\lfloor 2n/3 \rfloor}=0$.\\
                                                                      & $W_3: k_{3,4}=\pm \sqrt{-1/J_{\lfloor 2n/3 \rfloor}}\lambda^{-n/2}(1+ {\mathcal O}(\lambda^{-1/2}))$\\ \hline
\multicolumn{2}{|c|}{$g = 3q$, the Casimirs are $I_0=\prod_{j=0}^{n-1} a_{j,3};J_{\lfloor 2n/3 \rfloor}:=(-1)^n\prod_{j=0}^{n-1} a_{j,2};G_{\lfloor n/3 \rfloor}$.}\\ \hline
\multicolumn{2}{|c|}{$(\psi_{i,1}) \ge -D + 2O_2 + W_3 + i (W_3-O_2)$}\\
\multicolumn{2}{|c|}{$(\psi_{i,2}) \ge -D + O_2 + W_3 + i (W_3-O_2)$}\\
\multicolumn{2}{|c|}{$(\psi_{i,3}) \ge -D + W_1 + W_2 + W_3 + i (W_3-O_2)$}\\
\multicolumn{2}{|c|}{$(\psi_{i,4}) \ge -D + 3O_2 + i (W_3-O_2)$}\\ \hline
\end{tabular}
\end{center}
}

As an example, we show how to use these tables to find the genus of the spectral curve.
As before, we assume $n$ to be odd, $n=2q+1$. Recall  that for the short-diagonal  pentagram map $T_{2,1}$ 
 the genus is $g=3q$ for odd $n$, see  \cite{KS}.

\begin{proposition}\label{prop:Jacobians}
The spectral curves for the dented pentagram maps in $\CP^3$ generically have the genus $g=3q-1$ for $n$ odd and divisible by 3 and the genus $g=3q$ for $n$ odd and not divisible by 3.
\end{proposition}

\proof
Let us  compute the genus for the dented pentagram map $T_1$. As follows from the definition of the spectral curve $\Gamma$,
 it is a ramified 4-fold cover of $\CP^1$, since the $4\times 4$-matrix $\tilde{T}_{i,t}(\lambda)$ (or ${T}_{i,t}(\lambda)$) has 4 eigenvalues.
By the Riemann-Hurwitz formula the  Euler characteristic of $\Gamma$ is
$\chi(\Gamma)=4\chi(\CP^1)-\nu=8-\nu$, where $\nu $ is the ramification index of the covering.
In our setting, the index $\nu$ is equal to the sum of orders of the branch points at  $\lambda=0$ and
$\lambda=\infty$, plus the number  $\bar\nu$ of  branch points over $\lambda\not=0, \infty$, 
where  we   assume the latter points  to be all of order $1$ generically.
On the other hand, $\chi(\Gamma)=2-2g$, and once we know $\nu$ it allows us to find the genus
of the spectral curve $\Gamma$ from the formula $2-2g=8-\nu$.

The number $\bar\nu$ of branch points of $\Gamma$ on the $\lambda$-plane equals
the number of zeroes of the function $\partial_k R(\lambda,k)$
aside from the singular points $\lambda=0$ or $\infty$.
The function $\partial_k R(\lambda,k)$ is meromorphic on $\Gamma$,
therefore the number of its zeroes equals the number of its poles.
One can see that for any  $n=2q+1$  the function $\partial_k R(\lambda,k)$
has poles of total order $5n$ at $z=0$,
and it has zeroes of total order $2n$ at $z=\infty$.
Indeed, substitute the local  series for $k$ in $\lambda$ from the table to the expression for
$\partial_k R(\lambda,k)$.  (E.g., at $O_1$ one has $k={\mathcal O}(1)$.
The leading terms of $\partial_k R(\lambda,k)$ for the pole at $\lambda=0$ are
$4k^3, -3k^2G_0\lambda^{-q}, 2k J_0\lambda^{-n}, -I_0 \lambda^{-n}$. The last two terms, being  of order $\lambda^{-n}$, dominate and give the pole of order  $n=2q+1$.)
For $n$ odd and not divisible by 3, the corresponding orders of the poles and zeroes of $\partial_k R(\lambda,k)$ on  the curve
$\Gamma$ are summarized as follows:
\[
\begin{array}{||c|c||c|c||c|c||c|c||}
\hline
\text{ pole }  & \text{ order } & \text{ zero } & \text{ order }  \\ \hline
O_1     &  n            & W_1 & 0      \\ \hline
O_2     &  n            & W_{2} & 2n\\ \hline
O_{3} &  3n           & & \\ \hline
\end{array}
\]

Therefore, for such $n$, the total order of poles is $n+n+3n=5n$, while the total
order of zeroes is $0+2n=2n$.
Consequently, the number of zeroes of $\partial_k R(\lambda,k)$ at nonsingular points $\lambda\not=\{0,\infty\}$  is $\bar\nu=5n-2n=3n$, and so is
the total number of branch points of $\Gamma$ in the finite part of the $(\lambda,k)$ plane (generically, all of them have order $1$).
For $n$ odd and not divisible by 3 there is an additional branch point at $\lambda=0$ of order $1$ and another branch point at $\lambda=\infty$ of order $2$ (see the table for $T_1$).
Hence the ramification index is $\nu=\bar\nu+3=3n+3=6q+6$.
The identity $2-2g=8-\nu$ implies that $g=3q$.

For $n$ odd and divisible by 3, $n=6l+3$,  one has the same orders of poles $O_j$,  $W_1$ is of order zero, while each of the three zeros  $W_{2,3,4}$ is of order $4l+2$.  Then the total order of zeroes
is still $3(4l+2)=12l+6=2n$, and again $\bar\nu=5n-2n=3n$. However, there is no branch point at $\lambda=\infty$ and hence  the ramification index is $\nu=\bar\nu+1=3n+1=6q+4$. Thus for such $n$ we obtain from the identity $2-2g=8-\nu$ that $g=3q-1$.
\smallskip

Finally note that $T_1$ and $T^{-1}_2$ are conjugated to each other by means of the involution
$\alpha_{\mathbf 1}$, and hence  $T_1$ and $T_2$ have the same dimensions of invariant tori.
Their spectral curves are related by a change of coordinates furnished by this involution and have the same genus.
\proofend


\section{Dual dented  maps}\label{sect:dual}

\subsection{Properties of dual dented pentagram maps}

It turns out that the pentagram dynamics of $\hat{T}_m$ has the following simple description.
(We consider the geometric picture over $\R$.)

\begin{proposition}\label{prop:two_subs}
The dual pentagram map $\hat{T}_m$ in $\RP^d$ sends the vertex $v_k$ into the intersection of the
subspaces of dimensions $m$ and $d-m$ spanned by the vertices
$$
\hat{T}_m v_k =(v_{k+d-m-1}, ...,v_{k+d-1})\cap (v_{k+d}, ...,v_{k+2d-m})\,.
$$
\end{proposition}

\proof
As we discussed above, the point $\hat T_mv_k$ is defined by taking the intersection
of all but one consecutive hyperplanes:
$$
\hat T_mv_k=P_k\cap P_{k+1}\cap ...\cap   P_{k+d-m-1}\cap   P_{k+d-m+1}\cap ...\cap P_{k+d}.
$$
Note that this point $\hat T_mv_k$ can be described as the intersection of the subspace
$$
L_1^m =P_k\cap P_{k+1}\cap ...\cap  P_{k+d-m-1}
$$
of dimension $m$ and  the subspace
$$
L_2^{d-m} =P_{k+d-m+1}\cap ...\cap P_{k+d}
$$
of dimension $d-m$ in $\RP^d$. (Here the upper index stands for the dimension.)
Since each of the subspaces
$L_1$ and $L_2$ is the intersection of several consecutive hyperplanes $P_j$,
and each hyperplane $P_j$ is spanned by consecutive vertices, we see that
$L_1^m=(v_{k+d-m-1}, ...,v_{k+d-1})$ and  $L_2^{d-m}=(v_{k+d}, ...,v_{k+2d-m})$, as required.
\proofend

Consider the shift $Sh$ of vertex indices by $d-(m+1)$ to obtain the map
$$
\widehat{T}_m v_k := (\hat{T}_m \circ Sh)\, v_k=(v_{k}, ...,v_{k+m})\cap (v_{k+m+1}, ...,v_{k+d+1})\,,
$$
which we will study from now on.

\begin{figure}[hbtp]
\centering
\includegraphics[width=3.5in]{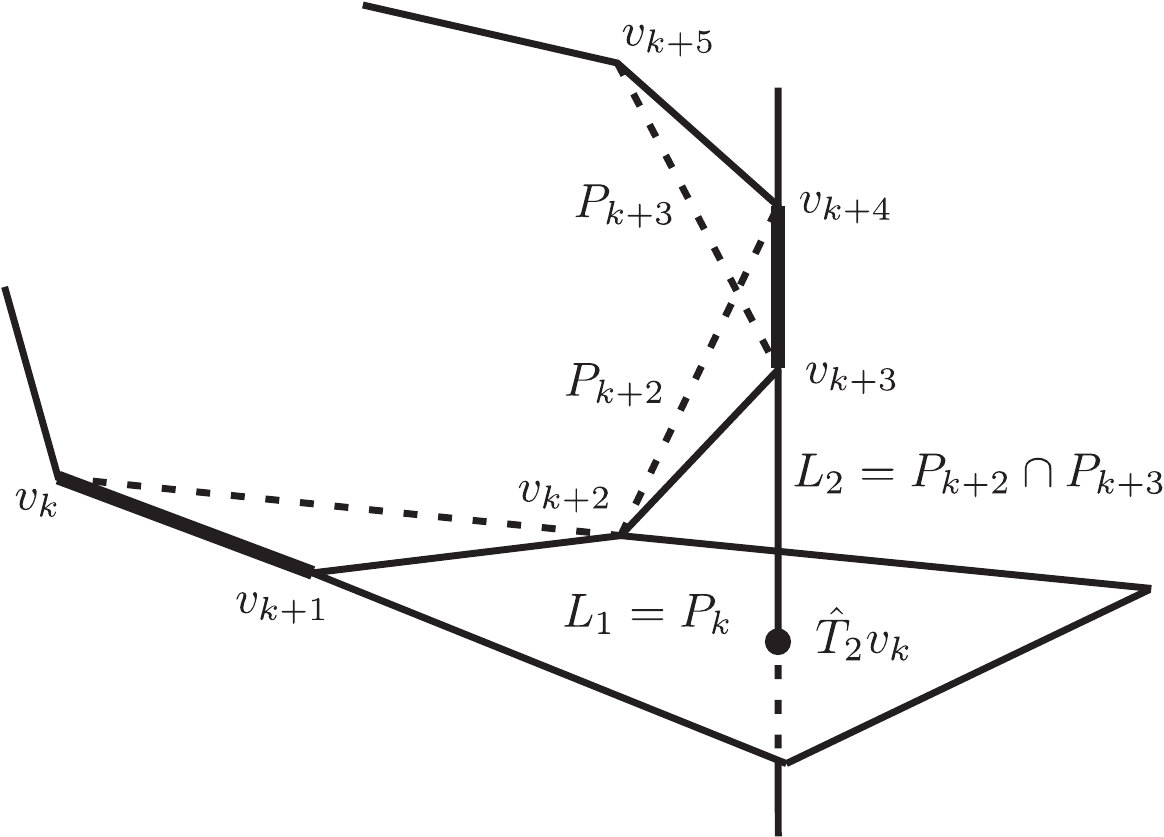}
\caption{\small The dual $\widehat T_2$ to the dented pentagram map $T_m$ for $m=2$ in $\RP^3$.}
\label{fig:dual-map}
\end{figure}

\begin{example}
{\rm
For $d=3$ and $m=2$ we have the dual pentagram map $\widehat{T}_2$ in $\RP^3$ defined via intersection of the 2-dimensional plane $L_1=(v_{k}, v_{k+1}, v_{k+2})$ and the line $L_2=(v_{k+3}, v_{k+4})$:
$$
\widehat{T}_2 v_k =(v_{k}, v_{k+1}, v_{k+2})\cap (v_{k+3}, v_{k+4})\,,
$$
see Figure \ref{fig:dual-map}. This map is dual to the dented pentagram map $T_m$ for $I=(1,2)$ and $J=(1,1)$.
}
\end{example}

\medskip

Let $V_k$ are the  lifts of the vertices $v_k$ of a twisted $n$-gon from $\RP^d$ to $\R^{d+1}$.
We assume that $n$ and $d+1$ are mutually prime and the conditions $\det(V_k,..., V_{k+d})=1$ with
$ V_{k+n}=MV_k$ for all $k\in \Z$ to provide the lift uniqueness.

\begin{proposition}\label{Tm-formula}
Given a twisted polygon $(v_k)$ in $\RP^d$ with coordinates $a_{k,j}$,  the image  $\widehat{T}_m V_k$ in $\R^{d+1}$ under the dual pentagram map is
proportional to the vector
$$
R_k=a_{k,m} V_{k+m} + a_{k,m-1} V_{k+m-1} +...+ a_{k,1} V_{k+1} +(-1)^{d} V_k
$$
for all $k\in \Z$.
\end{proposition}

\proof
Since
$$
\widehat{T}_m V_k \in (V_{k}, ...,V_{k+m})\cap (V_{k+m+1}, ...,V_{k+d+1})\,,
$$
the vector $W_k:=\widehat{T}_m V_k$ can be represented as a linear combination
of vectors from either of the groups:
$$
W_k=\mu_k V_{k}+ ...+\mu_{k+m}V_{k+m}=\nu_{k+m+1} V_{k+m+1}+ ...+\nu_{k+d+1}V_{k+d+1}.
$$
Normalize this vector by setting $\nu_{k+d+1}=1$. Now recall that
$$
V_{k+d+1} = a_{k,d} V_{k+d} + a_{k,d-1} V_{k+d-1} +...+ a_{k,1} V_{k+1} +(-1)^{d} V_k,
$$
for $k\in \Z$. Replacing $V_{k+d+1} $ by its expression via $V_k,...,V_{k+d}$ we obtain that
$\mu_k =(-1)^{d}$, $\mu_{k+1}=a_{k,1}, ... $, $\mu_{k+m}=a_{k,m}$. Thus  the vector
$$
R_k=a_{k,m} V_{k+m} + a_{k,m-1} V_{k+m-1} +...+ a_{k,1} V_{k+1} +(-1)^{d} V_k
$$
belongs to both the subspaces, and hence spans their intersection.
\proofend

Note that the image $\widehat{T}_m V_k$ under the dual map is
$W_k:=\widehat{T}_m V_k=\lambda_k R_k$, where the coefficients $\lambda_k$ are determined by the condition that $\det(W_k,..., W_{k+d})=1$ for all $k\in \Z$.


\subsection{Proof of the scale invariance}\label{sect:scale_proof}

In this section we prove scaling invariance in any dimension $d$ for any map
$\widehat T_m,\; 1 \le m \le d-1$ dual to the dented pentagram map $T_m$
on twisted $n$-gons in $\CP^d$, whose hyperplanes $P_k$ are defined by
taking consecutive vertices, but skipping the $m${\rm th} vertex.

\begin{theorem}{\bf (=\ref{thm:scaling}$\widehat~$)}\label{thm:scaling2}
The dual  dented pentagram map $\widehat T_m$ on twisted $n$-gons in $\CP^d$  is invariant with respect to the following scaling transformations:
$$
a_{k,1} \to s^{-1}a_{k,1} ,\; a_{k,2} \to s^{-2}a_{k,2} ,\;...\:, \; a_{k,m} \to s^{-m} a_{k,m},
$$
$$
 a_{k,m+1} \to s^{d-m}a_{k,m+1}, \;...\:, \;
a_{k,d} \to s a_{k,d}
$$
for all $s\in \C^*$.
\end{theorem}

\proof
The dual dented pentagram map is defined by $W_k:=\widehat{T}_m V_k=\lambda_k R_k$,
where the coefficients $\lambda_k$ are determined by the normalization condition: $\det(W_k,..., W_{k+d})=1$ for all $k\in \Z$.
The transformed coordinates are defined using the difference equation
$$
W_{k+d+1} = \hat{a}_{k,d} W_{k+d} + \hat{a}_{k,d-1} W_{j+d-1} +...+ \hat{a}_{k,1} W_{k+1}+(-1)^d W_k.
$$
The corresponding  coefficients $\hat{a}_{k,j}$ can be readily found using  Cramer's rule:
\begin{equation}\label{cram-rule}
\hat{a}_{k,j}
=\dfrac{\lambda_{k+d+1}}{\lambda_{k+j}}
\dfrac{\det(R_k,  R_{k+1}, ... ,R_{k+j-1},R_{k+d+1} , R_{k+j+1}, ... ,R_{k+d})}{\det(R_k , R_{k+1}, ... , R_{k+d})}
\end{equation}
The normalization condition reads as $\lambda_k \lambda_{k+1}...\lambda_{k+d}\det(R_k , R_{k+1}, ... , R_{k+d})=1$ for all $k\in \Z$.

To prove the theorem, it is sufficient to prove that the determinants in (\ref{cram-rule}) are homogenous in $s$ and find their degrees of homogeneity.

\begin{lemma}
The determinant
 $\det(R_k , R_{k+1}, ... , R_{k+d})$ has zero degree of homogeneity in $s$.
 The determinant in the numerator of formula (\ref{cram-rule}) has the same degree of homogeneity in $s$ as ${a}_{k,j}$.
\end{lemma}

The theorem immediately follows from this lemma since
even if $\lambda_k$ have some nonzero degree of homogeneity, it does not depend on $k$
anyway by the definition of scaling transformation, and it cancels out in the ratio. Hence the whole expression (\ref{cram-rule}) for
$\hat{a}_{k,j}$ transforms  just like  ${a}_{k,j}$, i.e., the dented pentagram map is invariant with respect to the scaling.

\prooflem
Proposition \ref{Tm-formula} implies that  the vector $R_k:=(V_k,V_{k+1},...,V_{k+d})\rr_k$ has an expansion
$$
\rr_k=((-1)^d,a_{k,1},...,a_{k,m},0,...,0)^t
$$
in the basis $(V_k,V_{k+1},...,V_{k+d})$, where $t$ stands for the transposed matrix.
Note that the vector $R_{k+1}$ has a similar expression
$\rr_k=((-1)^d,a_{k+1,1},...,a_{k+1,m},0,...,0)^t$
in the shifted basis
$(V_{k+1},V_{k+2},...,V_{k+d+1})$, but in the initial basis $(V_k,V_{k+1},...,V_{k+d})$
its expansion has the form $\rr_{k+1}=N_k\rr_k$ for the transformation matrix $N_k$
(see its definition in the proof of Theorem \ref{thm:lax_anyD}),
since the relation (\ref{eq:difference_anyD})
implies
$$
(V_{k+1},V_{k+2},...,V_{k+d+1})=(V_k,V_{k+1},...,V_{k+d})N_k\,.
$$

Note that formula (\ref{cram-rule}) is independent on the choice of the basis used and we expand vectors $R_k$ in the basis
$(V_{k+m+1},V_{k+m+2},...,V_{k+m+1+d})$. It turns out that the corresponding expansions
$\rr_k,...,\rr_{k+d+1}$ have particularly simple form in this basis, which is crucial for the proof.
We use hats, $\hat \rr_k,...\hat \rr_{k+d+1}$, when the vectors $R_k, ..., R_{k+d+1}$
are written in this new basis. Explicitly we obtain
\begin{align*}
\hat \rr_k &= (N_k N_{k+1} ... N_{k+m})^{-1} \rr_k=(-a_{k,m+1},-a_{k,m+2},...,-a_{k,d},1,0,...,0)^t, \\
\hat \rr_{k+1} &= (N_{k+1} N_{k+2} ... N_{k+m})^{-1} \rr_{k+1}=(0,-a_{k+1,m+1},-a_{k+1,m+2},...,-a_{k+1,d},1,0,...,0)^t,\\
&\ldots\\
\hat \rr_{k+m} &= N_{k+m}^{-1} \rr_{k+m}=(0,...,0,-a_{k+m,m+1},-a_{k+m,m+2},...,-a_{k+m,d},1)^t,\\
\hat \rr_{k+m+1} &= \rr_{k+m+1}=((-1)^d,a_{k+m+1,1},...,a_{k+m+1,m},0,...,0)^t,\\
\hat \rr_{k+m+2} &= N_{k+m+1} \rr_{k+m+2}=(0,(-1)^d,a_{k+m+2,1},...,a_{k+m+2,m},0,...,0)^t,\\
\hat \rr_{k+m+3} &= N_{k+m+1} N_{k+m+2} \rr_{k+m+3}=(0,0,(-1)^d,a_{k+m+3,1},...,a_{k+m+3,m},0,...,0)^t,\\
&\ldots\\
\hat \rr_{k+d+1} &= N_{k+m+1} N_{k+m+2}...N_{k+d} \rr_{k+d+1}=(0,...,0,(-1)^d,a_{k+d+1,1},...,a_{k+d+1,m})^t.
\end{align*}

Consider the matrix $\mathbf{M}=(\hat \rr_k,\hat \rr_{k+1},...,\hat \rr_{k+d+1})$ of size $(d+1) \times (d+2)$, which is essentially the matrix of the system of linear equations determining $\hat{a}_{k,j}$. All its entries are homogenous in $s$.
Also note that the determinant $\det(R_k , R_{k+1}, ... , R_{k+d})=\det(\hat \rr_k,\hat \rr_{k+1},...,\hat \rr_{k+d})$ is the minor formed by the first $d+1$ columns, while
the determinant in the numerator of formula (\ref{cram-rule}) is up to a sign the minor formed by crossing out the  $(j+1)$th column in $\mathbf{M}$.

For instance, for $d=6$ and $m=2$ this matrix has the form:
\[
\mathbf{M}=
\begin{pmatrix}
 -a_{k,3} &  0                 & 0                 &  1             &  0           &  0           &  0 &  0\\
-a_{k,4}  &  -a_{k+1,3}  &  0                & a_{k+3,1}  &  1            &  0          &  0 &  0\\
 -a_{k,5}  & -a_{k+1,4}  & -a_{k+2,3}  & a_{k+3,2} &  a_{k+4,1}  &  1         &  0 &  0\\
 -a_{k,6}  &  -a_{k+1,5} &  -a_{k+2,4} & 0              & a_{k+4,2} & a_{k+5,1} &  1 &  0\\
 1             &  -a_{k+1,6} &  -a_{k+2,5} & 0             & 0              & a_{k+5,2} & a_{k+6,1} &  1\\
0              &  1               & -a_{k+2,6} & 0               & 0             & 0              & a_{k+6,2} & a_{k+7,1}\\
0              &  0               & 1                & 0               & 0            &  0             & 0              & a_{k+7,2}
\end{pmatrix}.
\]

Let us form the corresponding matrix $\mathbf{D}$ of the same size representing the homogeneity degree of the entries of $\mathbf{M}$ given by the scaling transformations. One can assign  an arbitrary degree to a zero entry, and we do this in such a way that within each column the degrees  would change uniformly. Note that those degrees also change uniformly along all the rows as well, except for one simultaneous jump after the $(m+1)$th column.
In the above example one has
\[
\mathbf{D}=
\begin{pmatrix}
 4 &  5 & 6 &  0 &  1 &  2 &  3 &  4\\
 3 &  4 & 5 & -1 &  0 &  1 &  2 &  3\\
 2 &  3 & 4 & -2 & -1 &  0 &  1 &  2\\
 1 &  2 & 3 & -3 & -2 & -1 &  0 &  1\\
 0 &  1 & 2 & -4 & -3 & -2 & -1 &  0\\
-1 &  0 & 1 & -5 & -4 & -3 & -2 & -1\\
-2 & -1 & 0 & -6 & -5 & -4 & -3 & -2
\end{pmatrix}.
\]
Then the determinant of the minors obtained by crossing any of the columns are  homogeneous.
Indeed, elementary transformations on the matrix rows
by adding to one row another multiplied by a homogeneous in $s$ function preserve this table
of homogeneity degrees. On the other hand, producing the upper-triangular form by such transformations one can easily compute the homogeneity degree of  the corresponding minor.
For instance, the minor
$\det(\hat \rr_k,\hat \rr_{k+1},...,\hat \rr_{k+d})$  formed by the first $d+1$ columns
has zero degree. Indeed, we need to find the trace of the corresponding $(d+1)\times (d+1)$ degree matrix. It contains $m+1$ columns with the diagonal entries of degree $d-m$, as well as
$d-m$ columns with the diagonal entries of degree $-m-1$, i.e.,
the total degree is $(d-m)\cdot (m+1)+  (-m-1)\cdot(d-m)=0$. (In the example above it is
$4\cdot 3+(-3)\cdot 4=0$ on the diagonal for the first 7 columns.)
Similarly one finds the degree of any  $j$th minor of the matrix $\mathbf{M}$
for arbitrary $d$ and $m$ by calculating the difference of the degrees for the diagonal $(j,j)$-entry and the $(j,d+2)$-entry in the matrix $\mathbf{D}$.
\proofend

\begin{remark}
{\rm
The idea of using  Cramer's rule and a simple form of vectors $R_k$  was suggested 
in \cite{Beffa_scale} to prove the scale invariance
of short-diagonal  pentagram maps $T_{2,1}$. For the maps $\widehat{T}_m$ we 
employ  this approach along with  passing to the dual maps and the above ``retroactive" basis change. This 
choice of basis in the proof of Theorem \ref{thm:scaling2} also allows one to obtain explicit formulas for the pentagram map via the matrix  $\mathbf{M}$.
}
\end{remark}


\section{Continuous limit of general pentagram maps}\label{sect:cont_lim}

Consider the continuous limit of the dented pentagram maps on $n$-gons as $n\to\infty$.
In the limit   a generic twisted $n$-gon becomes a smooth non-degenerate quasi-periodic curve $\gamma(x)$ in $\RP^d$.
Its lift $G(x)$ to $\R^{d+1}$ is defined by the conditions that the components of the vector function
$G(x)=(G_1,...,G_{d+1})(x)$ provide the homogeneous coordinates for
$\gamma(x)=(G_1:...:G_{d+1})(x)$ in $\RP^d$ and $\det(G(x),G'(x),...,G^{(d)}(x))=1$ for all $x\in \R$.
Furthermore, $G(x+2\pi)=MG(x)$ for a given $M\in SL_{d+1}(\R)$. Then $G(x)$
satisfies the linear differential equation of order $d+1$:
$$
G^{(d+1)}+u_{d-1}(x)G^{(d-1)}+...+u_1(x)G'+u_0(x)G=0
$$
with periodic coefficients $u_i(x)$, which is a continuous limit of difference equation \eqref{eq:difference_anyD}. Here $'$ stands for $d/dx$.

Fix a small $\epsilon>0$ and  let $I$ be any $(d-1)$-tuple $I=(i_1,...,i_{d-1})$ of positive integers.
For the $I$-diagonal hyperplane
$$
P_k:=(v_k, v_{k+i_1}, v_{k+i_1+i_2},..., v_{k+i_1+...+i_{d-1}})
$$
its  continuous analog  is the hyperplane $P_\epsilon(x)$ passing through
$d$ points $\gamma(x),\gamma(x+i_1\epsilon),...,\gamma(x+(i_1+...+i_{d-1})\epsilon)$
of the curve $\gamma$. In what follows we are going to make a parameter shift in $x$
(equivalent to shift of indices) and define
$P_\epsilon(x):=(\gamma(x+k_0\epsilon),\gamma(x+k_1\epsilon),...,\gamma(x+k_{d-1}\epsilon))$, for any real $k_0<k_1<...<k_{d-1}$  such that $\sum_l k_l=0$.

Let $\ell_\epsilon (x)$ be the envelope curve for the family of hyperplanes $P_\epsilon(x)$  for a fixed $\epsilon$.
The envelope condition means that  $P_\epsilon(x)$ are the osculating hyperplanes of the curve $\ell_\epsilon (x)$, that is  the point $\ell_\epsilon (x)$ belongs to the hyperplane $P_\epsilon(x)$, while the vector-derivatives $\ell'_\epsilon (x),...,\ell^{(d-1)}_\epsilon (x)$ span this hyperplane
 for each $x$. It means that the lift of $\ell_\epsilon (x)$ to $L_\epsilon (x)$ in $\R^{d+1}$
satisfies the system of $d$ equations:
$$
\det ( G(x+k_0\epsilon), ..., G(x+k_{d-1}\epsilon),  L^{(j)}_\epsilon(x) )=0,\quad j=0,...,d-1.
$$

A continuous limit of the pentagram map  is defined as the evolution of the curve $\gamma$ in the direction of the envelope $\ell_\epsilon$, as $\epsilon$ changes.
Namely, one can show that the expansion of $L_\epsilon(x)$ has the form
$$
L_\epsilon(x)=G(x)+\epsilon^2 B(x)+{\mathcal O} (\epsilon^3)\,,
$$
where there is no term linear in $\epsilon$ due to the condition $\sum_l k_l=0$.
It satisfies the family of differential equations:
$$
L_\epsilon^{(d+1)}+u_{d-1,\ep}(x)L_\epsilon^{(d-1)}+...+u_{1,\ep}(x)L_\epsilon'+u_{0,\ep}(x)L_\epsilon=0, \text{ where } u_{j,0}(x)=u_{j}(x).
$$
Then the corresponding expansion of the coefficients $u_{j,\ep}(x)$ as $u_{j,\ep}(x)=u_{j}(x)+\ep^2w_j(x)+{\mathcal O}(\ep^3)$,
defines the continuous limit of the pentagram map as the system of evolution differential equations $du_j(x)/dt\, =w_j(x)$ for $j=0,...,d-1$. (This definition of limit assumes that we have the standard tuple $J={\mathbf 1}:=(1,...,1)$.)

\medskip

\begin{theorem}{\bf (Continuous limit)}\label{thm:cont}
The continuous limit of any generalized pentagram map $T_{I,J}$ for any $I=(i_1,...,i_{d-1})$ and $J={\mathbf 1}$
(and in particular, of any dented pentagram map $T_m$) in dimension $d$ defined by the system $du_j(x)/dt\, =w_j(x), \, j=0,...,d-1$ for $x\in S^1$  is the $(2, d+1)$-KdV flow of the Adler-Gelfand-Dickey  hierarchy on the circle.
\end{theorem}

\begin{remark}
{\rm
Recall that the $(n, d+1)$-KdV flow is defined on  linear differential operators
 $L= \partial^{d+1} + u_{d-1}(x) \partial^{d-1} + u_{d-2}(x) \partial^{d-2} + ...+ u_1(x) \partial + u_0(x)$  of order $d+1$ with periodic coefficients $u_j(x)$, where $\partial^{k}$ stands for $d^k/dx^k$.
One can define the fractional power
$L^{n/{d+1}}$ as a pseudo-differential operator for any positive integer $n$ and take
its pure differential part  $Q_n :=(L^{n/{d+1}})_+$. In particular, for $n=2$ one has $Q_2= \partial^2 + \dfrac{2}{d+1}u_{d-1}(x) $. Then the $(n, d+1)$-KdV equation is the  evolution equation on (the coefficients of) $L$ given by $dL/dt= [Q_n,L] $, see \cite{Adler}.

 For $d=2$ the  (2,3)-KdV equation is the classical Boussinesq equation  on the circle: $u_{tt}+2(u^2)_{xx}+u_{xxxx}=0$, which appears as the continuous limit of the 2D pentagram map \cite{OST99}.
}
\end{remark}

\proof
By expanding in the parameter $\epsilon$ one can show that $L_\epsilon(x)$ has the form
$L_\epsilon(x)=G(x)+\epsilon^2 C_{d,I}\left(\partial^2+ \dfrac{2}{d+1}u_{d-1}(x)  \right)G(x)+{\mathcal O} (\epsilon^3)$ as $ \ep\to 0$,  for a certain non-zero constant $C_{d,I}$, cf. Theorem~4.3 in \cite{KS}.
We obtain  the following evolution of the curve $G(x)$ given by the
$\epsilon^2 $-term of this expansion:
${dG}/{dt} = \left(\partial^2+ \dfrac{2}{d+1}u_{d-1}\right)G$, or which is the same,
${d}G/{dt} =Q_2G$.

We would like to find
the evolution of the operator $L$ tracing it.
For any $t$, the curve $G$ and the operator $L$ are related by the differential equation $LG=0$ of order $d+1$.
Consequently,
$d(LG)/dt=({d}L/{dt}) G +  L ({d}G/{dt})=0.$

Now note that if the operator $L$ satisfies the $(2,d+1)$-KdV equation
${d}L/{dt}=[Q_2, L]:=Q_2L-LQ_2,$
and $G$ satisfies ${d}G/{dt} =Q_2G$, we have the identity:
$$
\frac{dL}{dt} G +  L \frac{dG}{dt}=(Q_2L-LQ_2) G + L Q_2G= Q_2LG=0\,.
$$
In virtue of the uniqueness of the linear differential operator $L$ of order $d+1$
for a given fundamental set of solutions $G$, we obtain that indeed the evolution
of $L$ is described by the $(2,d+1)$-KdV equation.
\proofend


\section{Corrugated polygons and dented diagonals}\label{S:corrug}

\subsection{Pentagram maps for corrugated polygons}
In \cite{GSTV} pentagram maps were defined on spaces of corrugated polygons in $\RP^d$.
These maps turned out to be integrable, while the corresponding Poisson structures are related to many interesting structures on such polygons. Below we describe how one can view integrability in the corrugated case as a particular case of the dented maps.

Let $(v_k)$ be generic twisted $n$-gons in $\RP^d$ (here ``generic''   means that every $d+1$ consecutive vertices do not lie in a projective subspace). The space
of  equivalence classes of generic twisted $n$-gons in $\RP^d$ has dimension $nd$
and is denoted by ${\mathcal P}_n$.

\begin{definition}\label{def:2-corrugated}
A twisted polygon $(v_k)$ in $\RP^d$ is corrugated if for every $k\in \Z$ the vertices
$v_k, v_{k+1}, v_{k+d},$ and $v_{k+d+1}$ span a projective two-dimensional plane.
\end{definition}

The projective group preserves the space of corrugated polygons. Denote by
${\mathcal P}_n^{cor}\subset {\mathcal P}_n$ the space
of projective equivalence classes of generic corrugated $n$-gons.
One can show that such polygons form a submanifold of dimension $2n$ in the $nd$-dimensional space ${\mathcal P}_n$.

The consecutive  $d$-diagonals (the diagonal lines connecting
$v_k$ and $v_{k+d}$) of a corrugated polygon intersect pairwise, and the intersection points
form the vertices of a new corrugated polygon:   $T_{cor} v_k:=(v_k,v_{k+d})\cap(v_{k+1},v_{k+d+1})$.
This gives the definition of the pentagram map
on (classes of projectively equivalent) corrugated polygons $T_{cor}: {\mathcal P}_n^{cor}\to{\mathcal P}_n^{cor}$, see   \cite{GSTV}.
In 2D one has ${\mathcal P}_n^{cor}= {\mathcal P}_n$ and this gives the definition of the classical pentagram map on ${\mathcal P}_n$.

\begin{proposition}\label{prop:corr_well_def}{\rm \cite{GSTV}}
The pentagram map $T_{cor}$ is well defined on
${\mathcal P}_n^{cor}$, i.e., it sends a corrugated polygon to a corrugated one.
\end{proposition}

\proof
 The image of the pentagram map $T_{cor}$ is defined as the  intersection of the diagonals in the quadrilateral $(v_{k-1}, v_k, v_{k+d-1}, v_{k+d})$. Consider the diagonal $(v_k, v_{k+d})$.
 It contains  both vertices $T_{cor}v_{k-1}$ and $T_{cor}v_{k}$, as they are  
 intersections of this diagonal with the diagonals $(v_{k-1}, v_{k+d-1})$
and $(v_{k+1}, v_{k+d+1})$ respectively.
Similarly, both vertices $T_{cor}v_{k-d-1}$ and $T_{cor}v_{k-d}$ belong to the 
diagonal  $(v_{k-d}, v_{k})$.

Hence we obtain two pairs of new vertices $T_{cor}v_{k-d-1}, T_{cor}v_{k-d} $ and $T_{cor}v_{k-1}, T_{cor}v_{k}$ for each $k\in\Z$
 lying in one and the same 2D plane passing through old vertices $(v_{k-d}, v_{k}, v_{k+d})$. Note also that the indices of these new pairs  differ by $d$. Thus they satisfy the corrugated condition.
\proofend

\begin{theorem}\label{thm:restr_to_corr}
This pentagram map $T_{cor}: {\mathcal P}_n^{cor}\to {\mathcal P}_n^{cor}$ is a restriction of the dented pentagram map $T_m: {\mathcal P}_n\to {\mathcal P}_n$ for any $m=1,..., d-1$ from generic n-gons ${\mathcal P}_n$ in $\RP^d$ to corrugated ones ${\mathcal P}_n^{cor}$ (or differs from it by a shift in vertex indices).
\end{theorem}

In order to prove this theorem we first show that
the definition of a corrugated polygon in $\RP^d$ is equivalent to the following:

\begin{proposition}\label{prop:equiv_2-corrug}
Fix any $\ell=2, 3, ..., d-1$.
A generic twisted polygon $(v_k)$ is corrugated if and only if   the 2$\ell$ vertices
$v_{k-(\ell -1)}, ... , v_{k}$ and $v_{k+d-(\ell-1)}, ..., v_{k+d}$ span a projective $\ell$-space for every $k\in \Z$.
\end{proposition}

\proof
The case $\ell=2$ is the definition of a corrugated polygon. Denote the above
projective $\ell$-dimensional space by $Q^\ell_k=(v_{k-(\ell -1)}, ... , v_{k},v_{k+d-(\ell-1)}, ..., v_{k+d})$.

Then for any $\ell>2$ the intersections of the $\ell$-spaces
$Q^\ell_k$ and $Q^\ell_{k+1}$ is spanned by the vertices
$(v_{k-(\ell -2)}, ... , v_{k},v_{k+d-(\ell-2)}, ..., v_{k+d})$
and has dimension $\ell-1$, i.e., is the space $Q^{\ell-1}_k=Q^\ell_k\cap Q^\ell_{k+1}$.
This allows one to derive the condition on $(\ell-1)$-dimensional  spaces from the condition
on $\ell$-dimensional  spaces, and hence reduce everything to the case $\ell=2$.

Conversely, start with the  $(\ell-1)$-dimensional space
$Q^{\ell-1}_k$ and consider the space $Q^{\ell}_k$
containing $Q^{\ell-1}_k$, as well as the vertices $v_{k-(\ell -1)}$ and $v_{k+d-(\ell-1)}$.
We claim that after the addition of two extra vertices  the new space has dimension $\ell$, rather than $\ell+1$. Indeed,
the 4 vertices  $v_{k-(\ell -1)}, v_{k-(\ell -2)}, v_{k+d-(\ell-1)}, v_{k+d-(\ell-2)}$ lie in one and the same two-dimensional plane according to the above reduction. Thus adding  two vertices $v_{k-(\ell -1)}$ and $v_{k+d-(\ell-1)}$ to the space $Q^{\ell-1}_k$,
which already contains $v_{k-(\ell -2)}$ and $v_{k+d-(\ell-2)}$, boils down to adding one more
projective direction, because of the corrugated condition, and thus $Q^{\ell}_k$ has dimension $\ell$   for all $k\in \Z$.
\proofend

\proofthm
Now we take a generic twisted $n$-gon in $\RP^d$
and consider the dented $(d-1)$-dimensional diagonal $P_k$ corresponding to $m=1$ and
$I=(2,1,...,1)$, i.e., the hyperplane passing through the following $d$ vertices: $P_k=(v_k, v_{k+2}, v_{k+3},..., v_{k+d})$.

\begin{figure}[hbtp]
\centering
\includegraphics[width=6in]{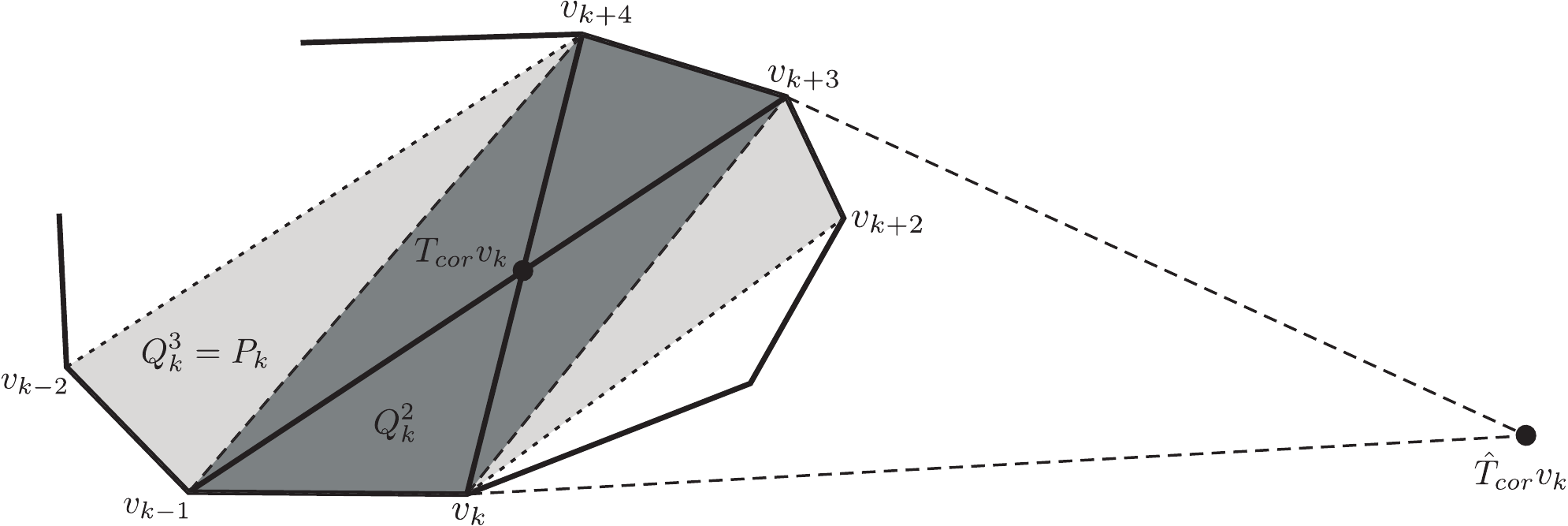}
\caption{\small The diagonal hyperplane $P_k$ coincides with the hyperplane $Q^{3}_k$ in $\RP^4$. Definitions of the corrugated pentagram map and its dual.}
\label{fig:corrugated}
\end{figure}

For a corrugated  $n$-gon in $\RP^d$, according to the proposition above, such a diagonal hyperplane will also pass through the vertices $v_{k-(\ell -1)}, ... , v_{k-1}$, i.e., it coincides with the space $Q^\ell_k$
for $\ell=d-1$:
$$
P_k=Q^{d-1}_k=(v_{k-(d-2)}, ... , v_{k},v_{k+2}, ..., v_{k+d})\,,
$$
see Figure \ref{fig:corrugated}.

Now the intersection of $d$ consecutive hyperplanes $P_k\cap P_{k+1}\cap...\cap P_{k+d-1}$, by 
the repeated use of the relation $Q^{\ell-1}_k=Q^\ell_k\cap Q^\ell_{k+1}$ for $\ell=d-1, d-2, ..., 3$,
reduces to the intersection of $Q^2_k\cap Q^2_{k+1}\cap Q^2_{k+2}$.
The latter is the  intersection of the diagonals in $Q^2_{k+1}$, i.e., $(v_{k+1},v_{k+d+1})\cap(v_{k+2},v_{k+d+2})=:T_{cor}v_{k+1}$. Thus the definition of the dented pentagram map $T_m$ for $m=1$ upon restriction to corrugated polygons reduces to the definition of the pentagram map $T$ on the latter (modulo shifts).

For any $m=1,...,d-1$ we consider the dented diagonal hyperplane
$$
P_{k-m+1}=(v_{k-m+1}, ..., v_k, v_{k+2}, v_{k+3},..., v_{k+d-m+1})\,.
$$
For corrugated  $n$-gons in $\RP^d$ this diagonal hyperplane  coincides with the space $Q^\ell_{k-m+1}$
for $\ell=d-1$ since it passes through all vertices from $v_{k-(d-2)}$ to $v_{k+d}$ with the exception of
$v_{k+1}$:
$$
P_{k-m+1}=Q^{d-1}_k=(v_{k-(d-2)}, ... , v_{k},v_{k+2}, ..., v_{k+d})\,.
$$
Thus the corresponding intersection of $d$ consecutive dented diagonal hyperplanes starting with
$P_{k-m+1}$ will differ only  by a shift of indices from the one for $m=1$.
\proofend

\begin{corollary}\label{cor:corr_coincide}
For dented pentagram maps $T_m$ with different values of $m$, their  restrictions
from generic to corrugated twisted polygons in $\RP^d$
coincide modulo an index shift.
\end{corollary}

Note that the inverse dented pentagram map $\widehat T_m$ upon restriction to corrugated polygons
also coincides with the inverse corrugated pentagram map $\widehat T_{cor}$. The latter is defined as follows:
for a corrugated polygon  $(v_k)$ in $\RP^d$
for every $k\in \Z$ consider the two-dimensional plane spanned by the vertices
$v_{k-1}, v_{k}, v_{k+d-1},$ and $v_{k+d}$. In this plane take the intersection of (the
continuations of) the sides on the polygon, i.e., lines  $(v_{k-1}, v_{k})$ and $(v_{k+d-1}, v_{k+d})$,
and set
$$
\widehat T_{cor}v_k:=(v_{k-1}, v_{k})\cap(v_{k+d-1}, v_{k+d})\,.
$$

\begin{corollary}
Continuous limit of the pentagram map $T_{cor}$
for corrugated polygons in $\RP^d$ is a restriction of the $(2, d+1)$-KdV equation.
\end{corollary}

The continuous limit for dented maps is found by means of the general procedure
described Section~\ref{sect:cont_lim}.
The restriction of the universal $(2, d+1)$-KdV system from generic
to corrugated curves might lead to other interesting equations on the submanifold.
(This phenomenon could be similar to  the KP hierarchy on  generic pseudo-differential
operators   $\partial+\sum_{j\ge 1} u_j(x)\partial^{-j}$, which when restricted
to operators of the form $\partial+\psi(x)\partial^{-1}\psi^*(x)$ gives the NLS equation, see \cite{Kr}.)

\medskip

\begin{remark}\label{rem:map_corr}
{\rm
One of applications of corrugated polygons is related to the fact that there is a natural map
from generic polygons in 2D to corrugated polygons in any dimension (see \cite{GSTV} and Remark
\ref{corrug-coord} below), which generically is a local diffeomorphism.
Furthermore, this map commutes with the pentagram map, i.e., it takes
deeper diagonals, which join vertices $v_i$ and $v_{i+p}$,
in 2D polygons to the intersecting diagonals of corrugated polygons in $\RP^p$.
This way one obtains a representation of the deeper diagonal pentagram
map $T_{p,1}$ in $\RP^2$ via the corrugated pentagram map in higher dimensions, see Figure \ref{fig:T31-2D}.
\begin{figure}[hbtp!]
\centering
\includegraphics[width=3.1in]{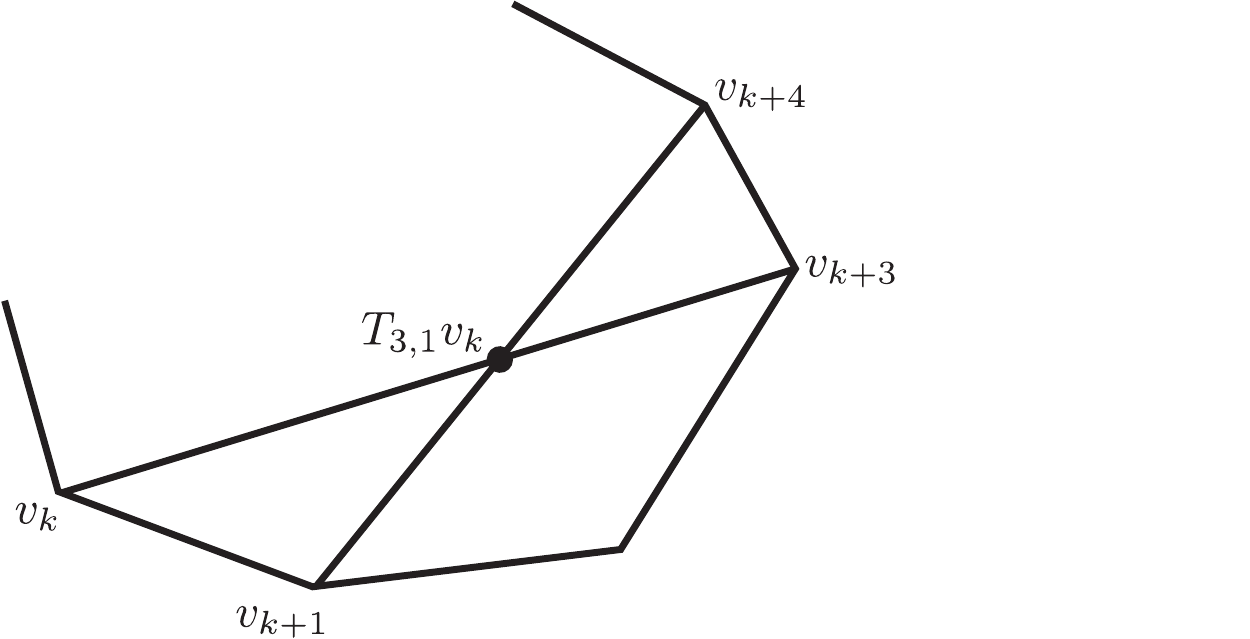}
\caption{\small Deeper pentagram map $T_{3,1}$ in 2D.} \label{fig:T31-2D}
\end{figure}

As a corollary one obtains that the  deeper diagonal pentagram map  $T_{p,1}$ in $\RP^2$ is
also an integrable system \cite{GSTV}.
Indeed, integrability of corrugated pentagram maps  implies  integrability of the pentagram map for deeper diagonals in 2D, since first integrals and other structures for the corrugated pentagram map in higher dimensions descends to those for the pentagram map on generic polygons in 2D
thanks to the equivariant local diffeomorphism between them.
 Explicit formulas of the invariants seems to be complicated because of a non-trivial relation between coordinates for polygons in $\RP^2$ and in $\RP^p$. 
}
\end{remark}


\subsection{Integrability for corrugated polygons}\label{corr3D}

Generally speaking, the algebraic-geometric integrability of the pentagram map on the space
${\mathcal P}_n$  (see Theorem \ref{thm:comparison} for 3D case), would not necessarily
imply the algebraic-geometric integrability for a subsystem, the pentagram map on
the subspace ${\mathcal P}_n^{cor}$ of corrugated polygons.

However, a Lax representation with a spectral parameter for corrugated polygons naturally follows from that for generic ones.
In this section, we perform its analysis in the 3D case (similarly to what has been done in Theorem \ref{thm:comparison}),
which implies the algebraic-geometric integrability for corrugated polygons in the 3D case.
It exhibits some interesting features: the dynamics on the Jacobian depends on whether $n$ is a multiple of $3$,
and if it is, it resembles a ``staircase'', but with shifts in $3$ different directions.
We also establish the equivalence of our Lax representation with that found in \cite{GSTV}.

For simplicity, we assume that $gcd(n,d+1)=1$ (see Remark \ref{diff-eq}).
In 3D case it just means that $n$ has to be odd.
Note that  this condition is technical, as it can be gotten rid of by using 
coordinates introduced in Section \ref{nonprimes}.\footnote{Another way to introduce the  coordinates is by means of the difference equation
$V_{j+d+1} = V_{j+d}+b_{j,d-1} V_{j+d-1}+...+b_{j,1} V_{j+1}+b_{j,0} V_j$,
used in \cite{GSTV}.}

\begin{remark}\label{corrug-coord}
{\rm
The coordinates on the space ${\mathcal P}_n^{cor}$ may be introduced using the same difference equation (\ref{eq:difference_anyD}) for $gcd(n,d+1)=1$. 
Since the corrugated condition means that the
vectors $V_j, V_{j+1}, V_{j+d}$ and $ V_{j+d+1}$ are linearly dependent for all $j \in \Z$, the subset
${\mathcal P}_n^{cor}$ of corrugated polygons  is singled out in
the space of generic twisted polygons  ${\mathcal P}_n$
by the relations $a_{j,l}=0, \; 2 \le l \le d-1$ in equation \eqref{eq:difference_anyD}, i.e., they are defined by the
equations
\begin{equation}\label{eq:2D_corr}
V_{j+d+1} = a_{j,d} V_{j+d} + a_{j,1} V_{j+1} +(-1)^{d} V_j,\quad j \in \Z\,.
\end{equation}

Furthermore, note that this relation also allows one to define a map $\psi$ from generic
twisted $n$-gons in $\RP^2$ to corrugated ones in $\RP^d$ for any dimension $d$,
see \cite{GSTV}. Indeed, consider  a lift of vertices $v_j \in \RP^2$ to vectors $V_j\in \R^3$
so that they satisfy the relations  \eqref{eq:2D_corr} for all $ j \in \Z$. Note that for $d\ge 3$ this is a nonstandard normalization of the lifts $V_j\in \R^3$, different from the one given in equation (\ref{eq:difference_anyD}) for $d=2$, since the vectors in the right-hand side are not consecutive.
Now by considering solutions $V_j\in \R^{d+1}$
of these linear relations \eqref{eq:2D_corr}
modulo the natural action of $SL_{d+1}(\R)$ we obtain a polygon in the projective space $\RP^d$ satisfying the corrugated condition. The constructed  map $\psi$ commutes with
the pentagram maps (since all operations are projectively invariant)
and is a local differeomorphism.
Observe that the subset ${\mathcal P}_n^{cor} \subset {\mathcal P}_n$  has dimension $2n$.
}
\end{remark}

Now we return to the consideration over $\C$.
The above restriction $gcd(n,d+1)=1$ allows one to define a Lax function in a straightforward way.
Here is an analogue of Theorem \ref{thm:comparison}:

\begin{theorem}
In dimension $3$ the subspace ${\mathcal P}_n^{cor} \subset {\mathcal P}_n$ is generically fibered into (Zariski open subsets of) tori of dimension
$g=n-3$ if $n=3l$, and $g=n-1$ otherwise.
\end{theorem}

\proof
The Lax function for the map $T_2$ restricted to the space ${\mathcal P}_n^{cor}$ is:
$$
 L_{j,t}^{-1}(\lambda) =
 \begin{pmatrix}
  0 & 0 & 0       & -1\\
  1 & 0 & 0       & a_{j,1}\\
  0 & 1 & 0       & 0\\
  0 & 0 & \lambda & a_{j,3}
 \end{pmatrix}.
$$
Now the spectral function has the form
$$
R(k,\lambda) = k^4 - \dfrac{k^3}{\lambda^{\lfloor n/3 \rfloor}} \left( \sum_{j=0}^{\lfloor n/3 \rfloor} G_j \lambda^j \right) +
 \dfrac{k^2}{\lambda^{\lfloor 2n/3 \rfloor}} \left( \sum_{j=0}^{N_0} J_j \lambda^j \right)
- \dfrac{k}{\lambda^n} \left( \sum_{j=0}^{\lfloor n/3 \rfloor} I_j \lambda^j \right) +\dfrac{1}{\lambda^n}
$$
where $N_0=\lfloor n/3 \rfloor-\lfloor gcd(n-1,3)/3 \rfloor$. One can show that 
$G_{\lfloor n/3 \rfloor}=\prod_{j=0}^{n-1} a_{j,1}$ and $I_0=\prod_{j=0}^{n-1} a_{j,3}$.
Below we present the summary  of relevant computations for the spectral functions, Casimirs, and
the Floquet-Bloch solutions, cf. Section \ref{sect:ag-in}.

\begin{center}
\begin{tabular}{|c|c|}
\hline
\multicolumn{2}{|c|}{$n=3l+1$; $n=3l+2$}\\ \hline
$\lambda=0$                                                & $\lambda=\infty$\\ \hline
$O_1: k_1 = 1/I_0 + {\mathcal O}(\lambda)$                         & $W_1: k_1=G_l(1+{\mathcal O}(\lambda^{-1}))$\\
$O_2: k_{2,3,4}= I_0^{1/3} \lambda^{-n/3} (1+{\mathcal O}(\lambda^{1/3}))$, & $W_2: k_{2,3,4}=G_l^{-1/3}\lambda^{-n/3}(1+ {\mathcal O}(\lambda^{-1/3}))$,\\ \hline
\multicolumn{2}{|c|}{$g = n-1$; there are $n+1$ first integrals; the Casimirs are $I_0,G_l$.}\\ \hline
\multicolumn{2}{|c|}{$(\psi_{i,1}) \ge -D + 2O_2 + W_2 +i(W_2-O_2)$}\\
\multicolumn{2}{|c|}{$(\psi_{i,2}) \ge -D + O_2 + W_1 + W_2 +i(W_2-O_2)$}\\
\multicolumn{2}{|c|}{$(\psi_{i,3}) \ge -D + W_1 + 2W_2 + i(W_2-O_2)$}\\
\multicolumn{2}{|c|}{$(\psi_{i,4}) \ge -D + 3O_2 + i(W_2-O_2)$}\\ \hline
\end{tabular}
\end{center}
\begin{center}
\begin{tabular}{|c|c|}
\hline
\multicolumn{2}{|c|}{$n=3l$}\\ \hline
$\lambda=0$                                                           & $\lambda=\infty$\\ \hline
$O_1: k_1 = 1/I_0 + {\mathcal O}(\lambda)$                            & $W_1: k_1=G_l(1+{\mathcal O}(\lambda^{-1}))$\\
$O_{2,3,4}: k_{2,3,4}= c_1 \lambda^{-l} (1+{\mathcal O}(\lambda))$,   & $W_{2,3,4}: k_{2,3,4}=c_2 \lambda^{-l}(1+ {\mathcal O}(\lambda^{-1}))$,\\
where $c_1^3-c_1^2 G_0+c_1 J_0-I_0=0.$                                & where $c_2^3 G_l-c_2^2 J_l+ c_2 I_l-1=0.$\\ \hline
\multicolumn{2}{|c|}{$g = n-3$; there are $n+3$ first integrals; the Casimirs are $I_0,G_0,J_0,I_l,J_l,G_l$.}\\ \hline
\multicolumn{2}{|c|}{$(\psi_{i,1}) \ge -D + O_2 + O_3 + W_4 +i_2(W_2-O_3) +i_1(W_3-O_2) +i_0(W_4-O_4)$}\\
\multicolumn{2}{|c|}{$(\psi_{i,2}) \ge -D + O_2 + W_1 + W_4 +i_2(W_2-O_2) +i_1(W_3-O_4) +i_0(W_4-O_3)$}\\
\multicolumn{2}{|c|}{$(\psi_{i,3}) \ge -D + W_1 + W_2 + W_4 +i_2(W_3-O_4) +i_1(W_4-O_3) +i_0(W_2-O_2)$}\\
\multicolumn{2}{|c|}{$(\psi_{i,4}) \ge -D + O_2 + O_3 + O_4 +i_2(W_4-O_4) +i_1(W_2-O_3) +i_0(W_3-O_2)$}\\
\multicolumn{2}{|c|}{where $i_2=\lfloor (i+2)/3 \rfloor, \quad i_1=\lfloor (i+1)/3 \rfloor, \quad i_0=\lfloor i/3 \rfloor$.}\\ \hline
\end{tabular}
\end{center}

The genus of spectral curves found above exhibits the dichotomy
$g=n-3$ or $g=n-1$ according to divisibility of $n$ by 3.
\proofend

\begin{remark}
{\rm
It is worth noting that  dimensions $g=n-3$ or $g=n-1$ of the Jacobians, and hence of the invariant tori, are consistent in the following sense:
\begin{itemize}
\item the sum of the genus of the spectral curve (which equals the dimension of its Jacobian) and the number of the first integrals equals $2n$, i.e.,
the dimension of the system;

\item the number of the first integrals minus the number of the Casimirs equals the genus of the curve. The latter also suggests that Krichever-Phong's universal formula provides a symplectic form for this system.
\end{itemize}

Also note that  Lax functions corresponding to the maps $T_1$ and $T_2$
restricted to the subspace ${\mathcal P}_n^{cor}$ lead to the same spectral curve in 3D,
as one can  check directly. This, in turn, is consistent with Corollary \ref{cor:corr_coincide}.
}
\end{remark}

\begin{proposition}
In any dimension the Lax function for the corrugated pentagram map $T_1$ in $\CP^d$ for  $gcd(n,d+1)=1$  is
\[
L_{j,t}(\lambda) =
\left(
\begin{array}{cccc|c}
0 & 0 & \cdots & 0    &(-1)^d\\ \cline{1-5}
\multicolumn{4}{c|}{\multirow{5}*{$D(\lambda)$}} & a_{j,1}\\
&&&& 0\\
&&&& \cdots\\
&&&& 0\\
&&&& a_{j,d}\\
\end{array}
\right)^{-1},
\]
with the diagonal $(d \times d)$-matrix $D(\lambda)={\rm diag}(1,\lambda, 1,...1)$. It is
equivalent to the one found in \cite{GSTV}.
\end{proposition}

\proof
The above Lax form follows from
Remark \ref{corrug-coord} and Theorem \ref{thm:lax_anyD}. To show the equivalence we define the gauge matrix as follows:
\[
g_j =
\left(
\begin{array}{cccc|c}
0 & 0 & \cdots & 0    &(-1)^d\\ \cline{1-5}
\multicolumn{4}{c|}{\multirow{4}*{$C_j$}} & 0\\
&&&& \cdots\\
&&&& 0\\
&&&& a_{j,d}\\
\end{array}
\right),
\]
where $C_j$ is the $(d \times d)$ diagonal matrix,
and its diagonal entries are equal to $(C_j)_{ll}=\prod_{k=0}^{d-l} a_{j-k,d},\; $ $1 \le l \le d$.
One can check that
\[
\tilde{L}_{j,t}(\lambda) = \dfrac{g_j^{-1} L_{j,t}^{-1} g_{j+1}}{a_{j+1,d}} =
\left(
\begin{array}{cccccc}
0       & 0 & 0 & \cdots & x_j & x_j+y_j\\
\lambda & 0 & 0 & \cdots & 0   & 0\\
0       & 1 & 0 & \cdots & 0   & 0\\
0       & 0 & 1 & \cdots & 0   & 0\\
\multicolumn{6}{c}{\cdots}\\
0       & 0 & 0 & \cdots & 1   & 1\\
\end{array}
\right),
\]
$$
\text{ where } x_j = \dfrac{a_{j,1}}{\prod_{l=0}^{d-1} a_{j-l,d}}, \quad y_j = \dfrac{1}{ \prod_{l=-1}^{d-1} a_{j-l,d} },
$$
which agrees with formula (10) in \cite{GSTV}.
\proofend

Note that the corresponding corrugated pentagram map has a cluster interpretation \cite{GSTV}
(see also \cite{Glick} for the 2D case).
On the other hand, it is a restriction of the dented pentagram map, which brings one to the following

\begin{problem}
Is it possible to realize the dented pentagram map $T_m$ on {\rm generic twisted polygons} in 
$\PP^d$  as a sequence of cluster transformations?
\end{problem}
We address this problem in a future publication.

\section{Applications:  integrability of pentagram maps for deeper dented diagonals}\label{sect:appl}

In this section we consider in detail more general dented pentagram maps.

\begin{definition}
{\rm
Fix an integer parameter $p\ge 2$ in addition to an integer parameter $m\in \{1,...,d-1\}$
and define the $(d-1)$-tuple  $I=I_{m}^{p}:=(1,...,1,p,1,...,1)$, where
the  value $p$ is situated at the $m$th place: $i_m=p$ and $i_\ell=1$ for $\ell\not=m$.
This choice of the tuple $I$ corresponds to the  diagonal plane $P_k$ which passes through
$m$ consecutive vertices $v_k, v_{k+1},...,v_{k+m-1}$, then skips $p-1$ vertices $v_{k+m}, ..., v_{k+m+p-2}$ (i.e., ``jumps to the next $p$th vertex")
and continues passing through the next $d-m$ consecutive vertices
$v_{k+m+p-1},...,v_{k+d+p-2}$:
$$
P_k:=(v_k, v_{k+1},...,v_{k+m-1},v_{k+m+p-1},v_{k+m+p},...,v_{k+d+p-2})\,.
$$
We call such a plane $P_k$ a {\it deep-dented diagonal (DDD) plane},
as the ``dent" now is of depth $p$, see Figure \ref{fig:DDD-plane}.
\begin{figure}[hbtp!]
\centering
\includegraphics[width=2.9in]{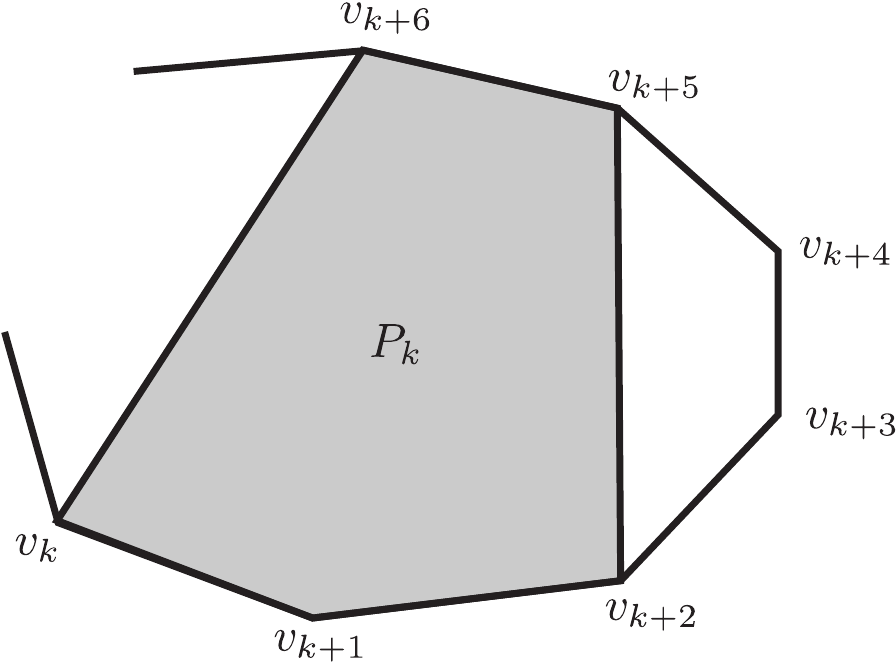}
\caption{\small The diagonal hyperplane for $I=(1,1,3,1)$ in $\RP^5$.} \label{fig:DDD-plane}
\end{figure}
Now we intersect $d$ consecutive planes $P_k$, to define the {\it deep-dented pentagram map} by
$$
T_m^p v_k:=P_{k}\cap P_{k+1}\cap ...\cap P_{k+d-1}\,.
$$
In other words, we keep the same definition of the $(d-1)$-tuple  $J={\mathbf 1}:=(1,1,...,1)$ as before: $T_m^p:=T_{I_{m}^{p},{\mathbf 1}}$.
}
\end{definition}

\begin{theorem}\label{thm:ddd}
The deep-dented  pentagram map for both twisted and closed polygons in any dimension
is a restriction of an integrable system to an invariant submanifold. Moreover, it admits a Lax representation with a spectral parameter.
\end{theorem}

To prove this theorem we introduce spaces of partially corrugated polygons, occupying intermediate positions between corrugated and generic ones.

\begin{definition}
{\rm
A twisted polygon $(v_j)$ in $\RP^d$ is {\it partially corrugated}
(or $(q,r;\ell)$-{\it corrugated}) if the diagonal subspaces $P_j$ spanned by two clusters
of $q$ and $r$ consecutive vertices $v_j$ with a gap of $(d-\ell)$ vertices between them
(i.e., $P_j=(v_j, v_{j+1},...,v_{j+q-1},v_{j+q+d-\ell},v_{j+q+d-\ell+1},...,v_{j+q+d-\ell+r-1})$, see Figure
\ref{fig:DDD-plane})  are subspaces of a fixed dimension $\ell\le q+r -2$ for all $j\in \Z$.
The inequality $\ell\le q+r -2$ shows that indeed these vertices are not in general position,
while $\ell= q+r -1$ corresponds to a generic twisted polygon.
We also assume that  $q\ge2$, $r\ge 2$, and $\ell\ge\max\{q,r\}$, so that the corrugated restriction would not be local, i.e., coming from one cluster of consecutive vertices, but would come from the interaction of the two clusters of those.
}
\end{definition}

Fix $n$ and denote the space of partially corrugated twisted $n$-gons  in $\RP^d$ (modulo projective equivalence) by
${\mathcal P}^{par}$.
Note that the corrugated condition in Definition \ref{def:2-corrugated} means $(2,2;2)$-corrugated
in this terminology.

\begin{proposition}\label{prop:equiv_partial-corrug}
The definition of a $(q,r;\ell)$-corrugated polygon in $\RP^d$ is equivalent to the
definition of a  $(q+1,r+1;\ell+1)$-corrugated polygon, i.e., one can add (respectively, delete) one extra vertex in each of the two clusters of vertices, as well as to increase (respectively, decrease) by one  the dimension of the subspace through them,
as long as $q,r\ge 2$,  $\ell \le q+r -2$, and $\ell\le d-2$.
\end{proposition}

\proof
The proof of this fact is completely analogous to the proof of Proposition \ref{prop:equiv_2-corrug} by adding  one vertex in each cluster.
\proofend

Define the {\it partially corrugated pentagram map} $T_{par}$ on the space ${\mathcal P}^{par}$:
to a partially corrugated twisted $n$-gon we associate a new one obtained by taking the intersections of $\ell+1$ consecutive diagonal subspaces $P_j$ of dimension $\ell$.

\begin{proposition}\label{prop:partial-pent}
i) The partially corrugated pentagram map is well defined:
by intersecting $\ell+1$ consecutive diagonal subspaces one generically
gets a point in the intersection.
ii) This map takes a partially corrugated polygon to a partially corrugated one.
\end{proposition}

\proof
Note that the gap of
$(d-\ell)$ vertices between clusters is narrowing by 1 vertex at each step as the dimension
$\ell$ increases by 1. Add maximal number of vertices, so that  obtain a hyperplane (of dimension $d-1$) passing through the clusters of $q$ and $r$ vertices with
a gap of one vertex between them. This is a dented hyperplane.
One can  see that intersections of 2, 3, ...  consecutive dented hyperplanes
gives exactly the planes of dimensions $d-2, d-3, ...$ obtained on the way while enlarging the clusters of vertices.
Then the intersection of $d$ consecutive dented hyperplanes is equivalent to the intersection of
$\ell+1$ consecutive  diagonal subspaces of dimension $\ell$ for partially corrugates polygons,
and generically is a point.

The fact that the image of a partially corrugated polygon is also  partially corrugated can be proved similarly to the standard corrugated case, cf. Proposition~\ref{prop:corr_well_def}.
We demonstrate the necessary changes in the following example. Consider the $(3,2;3)$-corrugated polygon in $\RP^d$ (here $q=3, r=2, \ell=3$), i.e., whose vertices
$(v_j, v_{j+1},v_{j+2},v_{j+d},v_{j+d+1})$ form a 3D
subspace $P_j$ in $\RP^d$ for all $j\in \Z$.
One can see that for the image polygon:
$a)$~three new vertices will be lying in the 2D plane obtained as the intersection  $B_{j+1}:=P_j\cap P_{j+1}=(v_{j+1},v_{j+2},v_{j+d+1})$ (since to get each of these three new vertices one needs to intersect this two planes with two more and the corresponding intersections will always be in lying in this plane);
$b)$~similarly, two new vertices will be lying on a certain line passing through the vertex $v_{j+d+1}$ (this line is the intersection of the 2-planes: $l_{j+d}:= B_{j+d}\cap B_{j+d+1}$).
Hence the obtained new 5 vertices belong to one and the same 3D plane spanned
by $B_{j+1}$ and  $l_{j+d}$ and hence satisfy the $(3,2;3)$-corrugated condition for all $j\in \Z$. The case of a general  partially corrugated condition is proved similarly.
\proofend

\begin{theorem}{\bf (=\ref{thm:ddd}$'$)}\label{thm:par}
The pentagram map on partially corrugated polygons in any dimension
is an integrable system: it admits a Lax representation with a spectral parameter.
\end{theorem}

\proofthms$\!\!\!\!${\bf \ref{thm:ddd} and \ref{thm:par}.}
Now suppose that we are given a generic polygon in $\RP^c$
and the pentagram map constructed with the help of a
deep-dented diagonals (of dimension $c-1$)  with the $(c-1)$-tuple of jumps $I=(1,...,1,p,1,...,1)$, which includes $m$ and $c-m$ consecutive vertices before and after the gap respectively.
Note that the corresponding gap between two clusters of points for such diagonals consists of $p-1$ vertices.
Associate to this polygon a partially $(q, r; \ell)$-corrugated polygon
in the higher-dimensional space
$\RP^d$ with clusters of $q=m+1$ and $ r=(c-m)+1$ vertices, the diagonal dimension is $\ell=c$,
and the space dimension is $d=c+p-2$.
Namely, in the partially corrugated polygon we add  one extra vertex  in each cluster,
increase the dimension of the diagonal plane by one as well (without the corrugated condition the diagonal dimension would increase by two after the addition of two extra vertices), while
the gap between two new clusters decreases by one: $(p-1)-1=p-2$. Then
the dimension $d$ is chosen so that the gap between two new clusters is $p-2=d-\ell$,
which implies that $d=\ell+p-2=c+p-2$.
(Example: for  deeper $p$-diagonals in $\RP^2$ one has $c=2, m=1, q=r=\ell=2$, and this way one obtains the space of corrugated polygons in $\RP^d$ for $d=p$.)
\medskip

Consider the map $\psi$ associating to a generic polygon in  $\RP^c$
a partially corrugated twisted polygon in $\RP^d$, where $d=c+p-2$.
(The map $\psi$ is defined similarly to the one for corrugated polygons in Remark
\ref{corrug-coord}.)
This map $\psi$  is a local diffeomorphism and commutes with the pentagram map:
the deep-dented pentagram map in $\RP^c$
is taken to the pentagram map $T_{par}$ on  partially corrugated  twisted polygons in $\RP^d$. In turn, the map $T_{par}$ is the restriction of the integrable dented pentagram map in $\RP^d$.
Thus the deep-dented pentagram map on polygons in $\RP^c$
is the restriction to an invariant submanifold of an integrable map on partially
corrugated  twisted polygons in $\RP^d$, and hence it is a subsystem of an integrable system.
The Lax form of the map $T_{par}$ can be obtained by restricting the Lax form for
dented maps from generic to partially corrugated polygons in $\RP^d$.
We present this Lax form  below.
\proofend

\begin{remark}\label{rem:coord_ddd}
{\rm Now we describe coordinates on the subspace of partially corrugated polygons and a Lax form of the corresponding pentagram map $T_{par}$ on them. 
Recall that on the space of generic twisted $n$-gons
$(v_j)$ in $\RP^d$ for  $gcd(n,d+1)=1$ there are coordinates $a_{j,k}$ for $0\le j \le n-1, 0\le k\le d-1$ defined by equation \eqref{eq:difference_anyD}:
$$
V_{j+d+1}=a_{j,d} V_{j+d} + ... + a_{j,1} V_{j+1} +(-1)^d V_j\,,
$$
where $V_j\in \R^{d+1}$ are lifts of vertices $v_j\in \RP^d$. One can see that the submanifold of
$(q,r;\ell)$-corrugated polygons in  $\RP^d$ without loss of generality can be assumed to have
the minimal number of vertices in clusters (see Proposition \ref{prop:equiv_partial-corrug}).
In other words, in this case there is such a positive integer $m$ that $q=m+1,  r=(\ell-m)+1$,
while the gap between clusters consists of $d-\ell$ vertices. Hence the corresponding
twisted polygons are described by linear dependence of
$q=m+1$ vertices $V_j, V_{j+1}, ..., V_{j+m}$ and $r$ vertices $V_{j+d+m-\ell+1},..., V_{j+d+1}$.
(Example: for $m=1, \ell=2$ implies a linear relation between $ V_j, V_{j+1}$ and $V_{j+d}, V_{j+d+1}$, which is the standard corrugated condition.)  This relation can be written as
$$
V_{j+d+1}=a_{j,d} V_{j+d} + ... +a_{j,d+m-\ell+1} V_{j+d+m-\ell+1}
+ a_{j,m} V_{j+m} +... +a_{j,1} V_{j+1} +(-1)^d V_j
$$
for all $j\in \Z$ by choosing an appropriate normalization of the lifts $V_j\in \R^{d+1}$.
Thus the set of partially corrugated polygons is obtained by imposing the condition
$a_{j,k}=0$ for $m+1\le k \le d+m-\ell$ and $0\le j \le n-1$ in the space
of generic twisted polygons given by equation \eqref{eq:difference_anyD}.
Note that the space of $(m+1, \ell-m+1 ;\ell)$-corrugated $n$-gons in  $\RP^d$
has dimension $n\ell$, while the space of generic twisted $n$-gons has dimension $nd$.
}
\end{remark}

In the complex setting, the Lax representation on such partially corrugated $n$-gons in  
$\CP^d$ or on generic
$n$-gons in  $\CP^c$ with deeper dented diagonals is described as follows.

\begin{theorem}
The  deep-dented pentagram map $T_m^p$ on generic
twisted and closed polygons in
$\CP^c$ and the  pentagram map $T_{par}$ on corresponding partially corrugated polygons
in $\CP^d$ with $d=c+p-2$ admits  the following Lax representation with a spectral parameter:
for $gcd(n,d+1)=1$ its Lax matrix is
\begin{equation}\label{eq:part_corr}
L_{j,t}(\lambda) =
\left(
\begin{array}{cccccc|c}
0 & 0 & \cdots & \cdots & 0  & 0  &(-1)^d\\ \cline{1-7}
\multicolumn{6}{c|}{\multirow{8}*{$D(\lambda)$}} & a_{j,1}\\
&&&&&& \cdots\\
&&&&&& a_{j,m}\\
&&&&&& 0\\
&&&&&& \cdots\\
&&&&&& 0\\
&&&&&&  a_{j,d+m-\ell+1}\\
&&&&&& \cdots\\
&&&&&& a_{j,d}\\
\end{array}
\right)^{-1},
\end{equation}
with the diagonal $(d \times d)$-matrix $D(\lambda)={\rm diag}(1,...,1,\lambda, 1,...1)$, where
$\lambda$ is situated at the $(m+1)${\rm th}  place, and an appropriate matrix $P_{j,t}(\lambda)$.
\end{theorem}

\proof This follows from the fact that the partially corrugated pentagram map is the restriction of the
dented map to the invariant subset of  partially corrugated polygons, so the Lax form is obtained by the corresponding restriction as well, cf. Theorem \ref{thm:lax_anyD}.
\proofend

Note that the jump tuple $I=(2,3)$ in $\PP^3$ corresponds to the 
first case which is neither a deep-dented  pentagram map, nor a short-diagonal one, and whose
integrability is unknown. It would be very interesting if the corresponding pentagram map turned out to be non-integrable. Some numerical evidence of non-integrability in that case is presented in \cite{KS14}.



\end{document}